\documentclass[11pt]{article}
\usepackage{latexsym}
\usepackage{graphics}
\usepackage{ifthen}
\usepackage{xspace,amssymb,latexsym,color}
\usepackage{amsmath}
\newcommand {\be}[1]{\begin{equation}\label{#1}}
\newcommand {\ee}{\end{equation}}
\newcommand {\bea}{\begin{eqnarray}}
\newcommand {\eea}{\end{eqnarray}}

\newcommand{\proof}{{\bf Proof : }}
\newcommand{\qed}{\hfill $\Box$}

\newcommand{\remarks}{{\bf Remarks : }}

\newcommand{\pr}{\mathbb{P}}

\newcommand{\bQ}{\bar Q}
\newcommand{\bT}{\bar T}
\newcommand{\bI}{\bar I}
\newcommand{\bA}{\bar A}
\newcommand{\bD}{\bar D}

\newcommand{\mb}[1]{\mbox{\boldmath $#1$}}
\newcommand{\mQ}{\mb{Q}}
\newcommand{\mT}{\mb{T}}

\newcommand{\mD}{\mb{D}}
\newcommand{\mZ}{\mb{Z}}
\newcommand{\mX}{\mb{X}}
\newcommand{\mS}{\mb{S}}

\newtheorem{theorem}{Theorem}
\newtheorem{lemma}[theorem]{Lemma}

\newtheorem{prop}{Proposition}
\newtheorem{coro}{Corollary}

\newtheorem{Defi}{Definition}

\title{\bf Instability in Stochastic and Fluid \\ Queueing Networks}

\author{
David Gamarnik
\thanks{IBM T.J. Watson Research Center, Yorktown Heights, NY 10598, USA.
        Email address: gamarnik@watson.ibm.com. }
\and
John J. Hasenbein
\thanks{Corresponding Author: The University of Texas at Austin,
Department of Mechanical Engineering,
1 University Station,
Austin, TX 78712-0292, USA.
        Email address: jhas@mail.utexas.edu.}
}

\begin{document}

\maketitle

\centerline{{\bf Keywords:} Rate stability, fluid models, large
deviations.}

\begin{abstract}
The fluid model has proven to be one of the most effective tools for the analysis of
stochastic queueing networks, specifically for the analysis of stability.
It is known that stability of a fluid model implies positive (Harris) recurrence (stability) of a
corresponding stochastic queueing network, and weak stability implies rate stability
of a corresponding stochastic network. These results have been established both for cases of specific
scheduling policies and for the class of all non-idling policies.

However, only partial converse results have been established and in certain cases converse statements
do not hold. In this paper we close one of the existing gaps. For the case of networks with two stations
we prove that if the fluid model is not
weakly stable under the class of all non-idling policies,
then a corresponding queueing network
is not rate stable under the class of all non-idling policies. We establish the result by
building a particular non-idling scheduling policy which makes the associated stochastic
process transient. An important corollary of our result is that the condition
$\rho^*\leq 1$, which was proven in \cite{daivan97} to be the exact condition for global
weak stability of the fluid model, is also the exact global rate stability condition
for an associated queueing network. Here $\rho^*$ is a certain computable
parameter of the network involving virtual station
and push start conditions.
\end{abstract}

\section{Introduction}\label{introduction}
In a series of papers, starting in the early 1990's, researchers established a
strong connection between the stability of a queueing network and the stability
of the corresponding fluid model. Initiated by Rybko and
Stolyar~\cite{rybsto92}
and generalized
by Dai~\cite{dai95b}, Stolyar~\cite{sto95}, and Chen~\cite{che95},
among others, it has been demonstrated that the stability of
a fluid model implies stability of a corresponding queueing network.
The stability results in the aforementioned papers were established
both for classes of policies, e.g.\ the set of non-idling policies,
and specific policies, e.g.\ First-In-First-Out (Dai~\cite{dai95b}
discusses both types of results).
The fluid model is a continuous,
deterministic analog of a discrete stochastic queueing network. It is
defined through
a set of equations which nominally take as parameters only the mean values
of the random variables
associated with the queueing network.
%In this paper when
%we refer to the fluid model it is understood to be the mean--value fluid model

Since the stability behavior of the fluid model is often significantly easier
to analyze than that of the stochastic model, the results above have led to
sweeping advances in understanding the stability of queueing networks via
the fluid model. A short list of such papers includes
\cite{bgt96,bra01,chezha94,daijen04,daivan97,kumkum96}.
However, a major element needed for a satisfactory theory of stability via
fluid models is a converse to the aforementioned stability results.
Specifically:
if the fluid model is \emph{not} stable in some sense, does this imply
instability
of the corresponding queueing network? Unfortunately, it turns out that
formulating
an appropriate converse is a delicate matter.
Partial converses which appear in the literature
refer both to the fluid model and the fluid limit model, which is the set
of weak
limits of the rescaled stochastic process.  Dai~\cite{dai96}
introduces the notion of a weakly unstable
fluid limit model. Roughly speaking, the fluid limit model is weakly
unstable if
there exists a uniform time at which all fluid limits which start at zero are
strictly positive. If the fluid limit model is weakly unstable,
Dai provides a concise proof showing
that in the stochastic network, the queue length process
diverges to infinity with probability one. This result provides a
partial converse to the stability results mentioned earlier.
Puhalskii and Rybko~\cite{puhryb00} use large deviations methods to prove
another partial
converse to the stability theorems.
Their result implies that if there exists an initial fluid model state for
which
all fluid trajectories with ``close'' initial states satisfy a uniform
rate of divergence condition, then the queueing process is not positive
Harris recurrent.
Under stronger conditions on the fluid trajectories they prove transience
of the queueing process. In two different papers, Meyn focuses on networks
which can be represented by countable state Markov chains.
In \cite{mey95}, Meyn uses martingale methods to show that
if all fluid limits
eventually diverge at some uniform rate, then the state process associated with
the queueing network is transient. Meyn~\cite{mey00} uses Markov chain
techniques to
prove another transience result. In that paper, if the fluid limits satisfy
a uniform homogeneity condition and a uniform lower bound for trajectories
starting from some open set, this again implies
that the state process for
the associated queueing network is transient. In each of the papers
\cite{mey95} and \cite{mey00}, Meyn
explains how the results can be extended to networks with more general
state spaces.

In all of the papers above which prove a converse result to the original
stability
theorems in \cite{che95,dai95b, sto95},
some uniform requirement
over a set of fluid trajectories, or more precisely a set of fluid
limits, (sometimes restricted to fluid limits starting from a particular
type of state) is needed for the result to be applicable.
Recall that
the original stability results of Dai~\cite{dai95b} and Chen~\cite{che95}
require that \emph{all} fluid trajectories are stable in some sense. Hence,
we use
the term ``partial converse'' above because there is some gap between the
stability and instability results.
%% The converses above do not describe what
%% the stability behavior of a network is if some fluid trajectories are stable
%% and some are unstable.
To close the gap between the stability and instability
results, one might consider analyzing directly the set of fluid limits.
However, this approach presents
certain difficulties, since the fluid limits are defined in a non-constructive
way, as weak limits of the underlying stochastic process. Moreover, it is shown
in Gamarnik~\cite{gar00} that computing fluid limits of a queueing
system is an algorithmically undecidable problem for a certain class of
scheduling
policies. In contrast, the fluid trajectories are defined by the set of
solutions
of a fluid model, a series of a fairly simple and reasonably tractable
differential equations.
It is this latter notion of fluid model which we use herein.

In Chen~\cite{che95}, it is shown that a multiclass queueing network
is globally rate stable if the corresponding fluid network is
globally weakly stable (see Section~\ref{subsubsection:StabilityFluid} in this paper for
definitions).
In this paper, we present a result which is
a full converse to Chen's stability result.
It is a full converse
in that for some networks, in particular two station networks, the result
implies that the stochastic network is globally rate stable if and only if
the corresponding fluid network is globally weakly stable. In particular,
this implies that if there is just one linearly divergent fluid trajectory,
then
the stochastic network is not rate stable under some non-idling policy.
Combining our main result
with the result of Dai and Vande Vate \cite{daivan97} we show that a
certain computable
condition of the form $\rho^*\leq 1$ is a necessary and sufficient
condition for rate stability
in networks with two stations. This is the first tight condition for
stability for such a general class
of networks. Our
proof uses a series of large deviations estimates to establish the
result and the only restriction in the stochastic network is that the
estimates are applicable to the primitive stochastic processes defining
the network. For a comprehensive discussion of various stability concepts
in fluid and queueing networks, and the $\rho^*\leq 1$ condition, we suggest
Dai~\cite{dai99}.

It should be noted that a strength of the transience
results in \cite{dai96,mey95,mey00,puhryb00}
is that they can be applied to networks under a class of
policies or just one particular policy (like FIFO or a static
buffer priority policy), whereas our result only applies to the
class of non-idling policies. In other words, the advantage
of the previous transience results is that they can be used
to determine if a given network is stable under a particular
scheduling policy. Our result can be used only to determine
if there exists one scheduling policy, within the class
of non-idling policies, which makes a network unstable. It should
be noted though that in general it is more difficult to apply
the previous results because of more stringent requirements
on the behavior of the fluid model trajectories.

One is naturally led to ask if our
result can be extended to apply to networks operating under a
particular policy rather than the class of all non-idling
policies. Unfortunately, a paper by Dai, Hasenbein, and
VandeVate~\cite{dhv00} essentially rules out the possibility of
obtaining a full converse which can also be applied to particular
policies. In that paper, it is shown that the stability of a
queueing network under a fixed static buffer priority policy
depends on more than just the mean value of the service and
interarrival times. Hence, no mean-value based fluid model can
sharply determine stability for the network considered, which
implies that no general stability converse can be formulated for a
network operating under an arbitrary, but specific policy.

Our paper is organized as follows, in Section~\ref{section:preliminaries} we
introduce stochastic and fluid multiclass networks and describe mathematical
preliminaries. In Section~\ref{section:main} we present the main results
of our paper and their implications. All of the proofs are presented in
Section~\ref{sec:proofs}.

\section{Preliminaries -- Model Description and Assumptions}\label{section:preliminaries}

We start by describing the model of interest - a multitype queueing
network. In the following subsection we describe a stochastic
multitype queueing network and in Subsection
\ref{subsection:FluidModel} we introduce a fluid queueing network.

\subsection{Stochastic multitype queueing networks}\label{subsection:stochQNET}

\subsubsection{Network description}
A open stochastic multitype queueing network is a network of $J$
stations $\sigma_1,\sigma_2,\ldots,\sigma_J$ each processing one
or multiple types of jobs. For each type $i=1,2,\ldots,I$ there is
an external stream of jobs arriving to the network. The
intervals between successive arrivals of jobs corresponding to
type $i$ are given by the i.i.d.\ sequence
$\mb{X}^i_1,\mb{X}^i_2,\ldots\,\mb{X}^i_k,\ldots.$
If $\mathbb{E}[\mX^i_1]$ exists, we
define $\lambda_i\equiv 1/\mathbb{E}[\mX^i_1]$ to be the arrival rate for type $i$.
More detailed assumptions about the stochastic processes  \{$\mb{X}^i_k,
k=1,2,\ldots\}$ are provided later. We denote by $\mb{A}_i(t)$ the
cumulative arrival process which counts the number of arrivals up
to time $t$. That is $\mb{A}_i(t)=\max\{k:\sum_{r\leq
k}\mb{X}^i_r\leq t\}$.

Each job of type $i$ has to be processed on a fixed ordered
sequence of stations
$\sigma(i,1),\sigma(i,2),\ldots,$ $\sigma(i,J_i)$, where each
$\sigma(i,l)$ is one of the stations $\sigma_1,\ldots,\sigma_J$.
We refer to $(i,1),(i,2),\ldots,(i,J_i)$ as stages corresponding
to the type $i$. We allow the repetition of stations, i.e.\
$\sigma(i,j')=\sigma(i,j'')$ for $j'\neq j''$, meaning some jobs need
to be processed on the same station multiple times (which is
common in some manufacturing environments). In particular $J_i$ could be
bigger than $J$. We slightly abuse the notation sometimes by using
$\sigma$ to also denote the set of classes which are served at
station $\sigma$.

Each station $\sigma=\sigma_j, j\leq J$ has one
server and, in particular, can work on only one job at a time.
Other jobs awaiting processing on $\sigma$ accumulate into queues.
Type $i$ jobs in the queue corresponding to stage $(i,j)$ will be
referred to as \emph{class $(i,j)$} jobs. Once a job of class
$(i,j)$ is processed, it is moved into the next queue $(i,j+1)$ at
the station $\sigma(i,j+1)$, or leaves the network if $j=J_i$. The
processing times for jobs of type $i$ at stage $j$ are random and
are given by the i.i.d.\ sequence
$\mb{S}^{i,j}_1,\mb{S}^{i,j}_2,\ldots,\mb{S}^{i,j}_k,\ldots.$
If $\mathbb{E}[\mS^{i,j}_1]$ exists, we define
$\mu_{i,j}\equiv 1/\mathbb{E}[\mS^{i,j}_1]$ to
be the service rate for jobs in class $(i,j)$.
Again, more detailed assumptions regarding the stochastic processes
$\{\mb{S}^{i,j}_k, k=1,2,\ldots\}$ are
provided later.

Let $d=\sum_i J_i$ denote the total number of classes in the
network. We denote by $\mb{Q}(t)=(\mb{Q}_{i,j}(t))\in
\mathbb{Z}^d_+$ the vector of queue lengths in our queueing
network at time $t\geq 0$. In order to completely specify the
stochastic dynamics of $\mb{Q}(t)$ we need to specify  the
vector of initial queue lengths $q=\mb{Q}(0)$ and the
\emph{scheduling policy} ${\cal U}$ which gives gives the protocol
at each station $\sigma$ for resolving the contention for service,
when several jobs are competing for the same station. Some common
policies include the First-In-First-Out (FIFO) policy which gives
priority to jobs which arrived earlier to the station,
Last-In-First-Out (LIFO) defined analogously, Global-FIFO (GFIFO)
which gives priority to jobs which arrived earlier into the entire
network (based on time stamps of a job's arrival to class
$(i,1)$), and static buffer priority policies which are based on a
ranking of classes in each station and give priority to jobs with
the higher ranking, etc. All of these policies are examples of
\emph{non-idling} policies, which are defined as policies that
require each station $\sigma$ to work at full capacity as long as
there are any jobs waiting to be processed by $\sigma$.

Throughout the paper
we will be only considering head-of-the-line (HOL) type non-idling scheduling policies.
Under a HOL policy, at most one job of each class at a given station can
receive service at a given time. Furthermore under the HOL assumption, jobs
are served in FIFO order within a given class.
FIFO, Global-FIFO and static buffer priority are examples
of HOL policies.
Adopting the HOL assumption in this paper
is really not a restriction since the main goal of this paper is to construct
an unstable (in a sense to be defined) non-idling scheduling policy.
Indeed we construct an unstable
policy which happens to be of HOL type.
In addition to being HOL, the policy we use to prove the main result
is \emph{preemptive resume}. Under such a policy if the processing
of a job of class $i$ is interrupted to serve a job from another
class $j$, then the class $i$ job is ejected from service and placed
at the head of the line for processing at a later time.
When class $i$ is again chosen for service, the remaining processing
time for the ejected job is the same as it was at the moment it was
ejected.

For each
$(q,z_1,z_2)\in\mathbb{Z}^d_+\times\Re_+^{I+d}$ we say that the
state of the stochastic process at time $t$ is $(q,z_1,z_2)$ if at
time $t$ the vector of queue lengths $\mb{Q}(t)$ is $q$, the
vector of residual interarrival times is $z_1$ (hence the
dimension $I$ for this component of the state) and the vector of
residual service times is $z_2$. For many scheduling policies, including the policy
constructed in this paper, the state space
$\mathbb{Z}^d_+\times\Re_+^{I+d}$ is adequate to describe the underlying stochastic
process of the network.

For each class $(i,j)$, let $\mb{T}_{i,j}(t)$ denote the total
amount of time station $\sigma(i,j)$ spent processing class
$(i,j)$ jobs during the time interval $[0,t]$. Let
$\mb{D}_{i,j}(t)$ denote the cumulative departure process for
class $(i,j)$ jobs, that is $\mb{D}_{i,j}(t)$ is the number of
class $(i,j)$ jobs that station $\sigma_{i,j}$ processed during
the time interval $[0,t]$. For each
station $\sigma$, let
$\mb{Q}_{\sigma}(t)=\sum_{(i,j)\in\sigma}\mb{Q}_{i,j}(t)$ and let
$\mb{T}_{\sigma}(t)=\sum_{(i,j)\in\sigma}\mb{T}_{i,j}(t)$. The
following relations follow immediately from the definitions. For
all $1\leq i\leq I, 2\leq j\leq J_i$ and $t\geq 0$
\begin{align}
\mb{Q}_{i,1}(t)&=\mb{Q}_{i,1}(0)+\mb{A}_i(t)-\mb{D}_{i,1}(t),  \label{eq:stQueueStage1} \\
\mb{Q}_{i,j}(t)&=\mb{Q}_{i,j}(0)+\mb{D}_{i,j-1}(t)-\mb{D}_{i,j}(t),  \label{eq:stQueueStagej} \\
\mb{D}_{i,j}(t)&=\max\{k:\sum_{r\leq k}\mb{S}^{i,j}_r\leq
\mb{T}_{i,j}(t)\}.  \label{eq:stTimeDeparture}
\end{align}

Also for every $0 \leq t_1\leq t_2$ and every station $\sigma$
\be{eq:stFeasibility}
\sum_{(i,j)\in\sigma}(\mb{T}_{i,j}(t_2)-\mb{T}_{i,j}(t_1))\leq
t_2-t_1. \ee Finally, if the scheduling policy ${\cal U}$ is
non-idling, then for every $0 \leq t_1\leq t_2$ and every station
$\sigma$, if $\mb{Q}_{\sigma}(t)>0$ for all $t_1\leq t\leq t_2$
then $\mb{T}_{\sigma}(t_2)-\mb{T}_{\sigma}(t_1)=t_2-t_1$. In other
words, if the total queue in station $\sigma$ was always positive
during the time interval $[t_1,t_2]$, then the station was always
working on jobs full time during this interval.

Let
\be{eq:mulambdamax}
\lambda_{\max}=\max_{i}\{\lambda_i,\lambda_i^{-1}\},
\ee
\be{eq:mumax}
\mu_{\max}=\max_{i,j}\{\mu_{i,j},\mu_{i,j}^{-1}\},
\ee
\be{eq:Jmax}
J_{\max}=\max_i\{J_i\}.
\ee
For technical purposes we introduce $C$ -- a very large constant which
exceeds all the parameters of the network. Specifically,
\be{eq:C}
C>13(\lambda_{\max}+\mu_{\max})^2IJ_{\max}^3.
\ee
For any station $\sigma$, let $|\sigma|$ denote the number of classes in the set $\sigma$.
For any vector $q\in \Re^d$ we let
$||q||=\sum_{1\leq i\leq d}|q_i|$ denote the $L_1$ norm. For any non-decreasing non-negative function
$f(t)$ and any $t_1\leq t_2$ we let $f(t_1,t_2)$ denote $f(t_2)-f(t_1)$.

\subsubsection{Stochastic Assumptions}
\label{subsubsection:LargeDeviations}

Below, we introduce some basic assumptions on the sequences
of random variables which represent the primitive data in
our stochastic networks, and an assumption on the behavior of
the network process itself.

\begin{Defi}\label{definition:LD}
Consider a sequence of i.i.d.\ non-negative random variables
$\mZ_1,\mZ_2,\ldots,\mZ_n, \ldots$ with
$\mathbb{E}[\mZ_1]\equiv\alpha < \infty$.
Such a sequence satisfies
large deviations (LD) bounds if for every $\epsilon>0$ there exist
constants $L=L(\epsilon),V=V(\epsilon)>0$ such that for any $z>0$
\be{eq:LDBtime} \mathbb{P}\Bigg(\Big|\sum_{1\leq i\leq
n}\mb{Z}_i-z-\alpha n \Big| \ge \epsilon n \; \Big| \; \mZ_1\geq z\Bigg)\leq Ve^{-Ln},
\ee for all $n\geq 1$, and the counting process
$\mb{N}(t)\equiv\max\{n:\mb{Z}_1+\cdots+\mb{Z}_n\leq t\}$
satisfies
\be{eq:LDBrate} \mathbb{P}\Bigg( \Big|\mb{N}(t+z)-{t\over
\alpha}\Big| \ge \epsilon t \; \Big| \; \mZ_1\geq z\Bigg)\leq Ve^{-Lt},
\ee
for all $t\geq 0$.
\end{Defi}

It is important that the constants $L,V$ in the definition above do
not depend on $z>0$. This uniformity will become useful when we
analyze arrival and service processes with the presence of some
residual interarrival and service times. For simplicity we assume
common constants $L=L(\epsilon),V=V(\epsilon)$ instead of
individual constants corresponding to indices $i,j$.
Our main stochastic
assumptions are as follows.

\begin{itemize}

\item {\bf Assumption A:} The sequences \{$\mX^i_n, n=1,2,\ldots \}$ and
$\{ \mS^{i,j}_n, n=1,2,\ldots \}$ are i.i.d.\ for every
$1\leq i\leq I$ and $1\leq j\leq J_i$.

\item{\bf Assumption B:}
For each $i,j$
the large deviation bound
(\ref{eq:LDBtime}) holds for the sequences
\{$\mX^i_n, n=1,2,\ldots \}$ and
$\{ \mS^{i,j}_n, n=1,2,\ldots \}$ and
(\ref{eq:LDBrate}) holds for the associated renewal processes.

\item {\bf Assumption C:} For every state $(q,z_1,z_2)$, every $n>0$ we have
$\tau\equiv\inf\{t:||\mb{Q}(t)||\geq n \; | \;\mb{Q}(0)=(q,z_1,z_2)\}<\infty$
with probability one, under any scheduling policy.

\end{itemize}

One way to verify assumption B is via the following sufficient
condition for the large deviations bounds to hold for
i.i.d.\ sequences.

\begin{lemma}\label{lemma:LDiid}
Suppose $\mZ_1,\mZ_2,\ldots,\mZ_m,\ldots$, is a non-negative
i.i.d.\ sequence with $\mathbb{E}[\mZ_1]=\alpha$ such
that there exists a function $F(\theta),\theta \geq 0$ taking values in
$\Re_+\cup\{\infty\}$
which is finite on some interval $[0,\theta_0]$ and which satisfies
\be{eq:exptail}
\sup_{z\geq 0}\mathbb{E}[e^{\theta (Z_1-z)}\; | \; Z_1\geq z]\leq
F(\theta),
\ee
for every $\theta\geq 0$.
Then this sequence satisfies the LD bounds (\ref{eq:LDBtime}) and
(\ref{eq:LDBrate}).
\end{lemma}

The proof of Lemma~\ref{lemma:LDiid}  is provided in the Appendix.
It is simple to check that condition (\ref{eq:exptail}) is
satisfied by many distributions including the exponential, Erlang,
and any distribution with bounded support. Note that, by setting $z=0$,
condition (\ref{eq:exptail}) implies that the distribution of
$Z_1$ has a moment generating function for $\theta \in
[0,\theta_0]$.

Assumption C is intentionally broad, in that it does not involve
the stochastic primitives directly. The assumption holds for a
wide range of distributions, given assumptions A and B.
For example, if at least one of
the service time distributions has unbounded support, then assumption
C holds. However, assumption C holds under even weaker conditions.

We adopt assumptions A, B, and C for the remainder of the paper.
Whenever we talk about the probability
$\mathbb{P}\{\cdot\}$ of any event, the probability is understood
with respect to the stochastic processes $\{ \mX^i_n,\mS^{i,j}_n \}$. If
the vector of initial queues $\mQ(0)$ is a random vector itself,
then the probability is also with respect to the probability
distribution of $\mQ(0)$.

\subsubsection{Stability and rate stability}\label{subsubsection:stability}
One of the main features one desires to have in a multitype
queueing network is stability.  Various equivalent
definitions of stability have been used in the literature, among
which positive Harris recurrence is one of the most commonly used
definitions. Under the condition that the interarrival times
$\{X_k^i \}$ are
unbounded and spread out (see \cite{dai95b}), then positive
Harris recurrence is defined as follows.

\begin{Defi}\label{Defi:Stable}
A multitype queueing network operating under a scheduling policy
${\cal U}$ is defined to be Harris recurrent if there exists $b>0$
such that for any initial vector of queue lengths $\mQ(0)$, the
time $\tau=\inf\{t:||Q(t)||\leq b\}$ is finite with probability
one. The network is defined to be positive Harris recurrent or
stable, if in addition $\mathbb{E}[\tau]<\infty$, where the
expectation is conditioned on the initial vector of queue lengths $\mQ(0)$. The
network is defined to be globally stable if it is stable for every
non-idling scheduling policy ${\cal U}$.
\end{Defi}
The positive Harris recurrence property, under some additional
technical assumptions, implies the existence of a unique
stationary distribution for the queue length process $\mQ(t)$.

A somewhat weaker definition of stability is rate stability. This
is the form of stability we are primarily concerned with in
this paper.

\begin{Defi}\label{Defi:RateStable}
A multitype queueing network operating under a scheduling policy
${\cal U}$ is defined to be rate stable if for every type $i$,
$\lim_{t\rightarrow\infty}{\mD_{i,J_i}(t)\over t}=\lambda_i$, a.s.
The network is defined to be globally rate stable if it is rate
stable for every non-idling scheduling policy ${\cal U}$.
\end{Defi}

In words, rate stability means with probability one the arrival rate is equal
to the departure rate. From (\ref{eq:stQueueStage1}) and
(\ref{eq:stQueueStagej}), rate stability implies
$\lim_{t\rightarrow\infty}{\mD_{i,j}(t)\over t}=\lambda_i$ and
$\lim_{t\rightarrow\infty}{\mQ_{i,j}(t)\over t}=0$ a.s.\ for all
$i,j$. In other words, for a rate stable system,
even if the total queue length $||\mQ(t)||$ diverges as $t$ goes
to infinity, it grows at most at a sub-linear rate a.s.

\subsection{Fluid model}\label{subsection:FluidModel}
\subsubsection{Fluid equations}\label{subsubsection:FluidEquations}
Fluid models are continuous deterministic counterparts of
stochastic queueing networks, intended to capture the most
essential dynamic properties of the queue length process.
The term fluid model is sometimes used interchangeably with
the terms ``fluid limits'' and ``functional law
of large numbers.''
For many types of
queueing networks (see, e.g.\ \cite{dai95b},
%\cite{rs}, \cite{bramsoncounter}, \cite{cof_stolyar}, \cite{dantzer_robert3}
\cite{rybsto92}, \cite{bra98c}, \cite{cofsto01}, \cite{danrob02})
it has been established that
the rescaled queue length process $\mQ(nt)/n$ for a large
scaling parameter $n$ converges weakly to a certain continuous
deterministic process, satisfying a series of functional
equations, which we describe below. To avoid confusion we
define the \emph{fluid limit model} to be the set of weak limits of
 $\mQ(nt)/n$ as $n\rightarrow\infty$, and we define the fluid
model to be the set of solutions of the system of equations below (formal definition follows).
Then the set of fluid limits is a subset of the set of solutions to the fluid model.

Given a multitype queueing network with arrival rates $\lambda_i$
and service rates $\mu_{i,j}$, the corresponding HOL fluid model (or fluid network) is
defined by the following system of equations and inequalities
with time dependent variables  $\bQ_{i,j}(t), \bA_i(t),
\bD_{i,j}(t),\bT_{i,j}(t), t\geq 0$. We first provide the system of
equations and inequalities, and immediately after we give a
physical explanation for each of these equations. For every
$i=1,\ldots,I, j=1,2,\ldots,J_i$,
$\sigma=\sigma_1,\ldots,\sigma_J$, $t\geq 0$ and $0\leq t_1\leq t_2$
\begin{align}
&\bQ_{i,1}(t)=\bQ_{i,1}(0)+\bA_i(t)-\bD_{i,1}(t), \label{eq:flStage1} \\
&\bQ_{i,j}(t)=\bQ_{i,j}(0)+\bD_{i,j-1}(t)-\bD_{i,j}(t), \label{eq:flStagej} \\
&\bA_i(t)=\lambda_it, \label{eq:flA=lambda t} \\
&\bD_{i,j}(t)=\mu_{i,j}\bT_{i,j}(t), \label{eq:flD=muT} \\
&\sum_{(i,j)\in\sigma}(\bT_{i,j}(t_2)-\bT_{i,j}(t_1))\leq t_2-t_1, \label{eq:flFesibility} \\
&\bQ_{i,j}(t),\bA_{i}(t),\bD_{i,j}(t),\bT_{i,j}(t)\in \Re_+. \label{eq:flNonnegativity}
\end{align}
In addition, for all $i,j, \bT_{i,j}(t)$ is  a non-decreasing
function of $t$ and $\bA_{i,j}(0),\bT_{i,j}(0)=0.$

The value of $\bQ_{i,j}(t)$ represents the total amount of fluid
present in buffer $(i,j)$ at time $t$. We also refer to it
as class $(i,j)$ fluid. $\bA_i(t)$ represents the total amount of
fluid corresponding to type $i$, that arrived externally during
the time interval $[0,t]$. The fluid arrival process is assumed to
be linear with rate $\lambda_i$, hence (\ref{eq:flA=lambda t}).
$\bD_{i,j}(t)$ is the amount of class $(i,j)$ fluid that
was processed by station $\sigma(i,j)$ during $[0,t]$.
$\bT_{i,j}(t)$ represents
the portion of the time interval $[0,t]$ that station
$\sigma(i,j)$ spent processing class $(i,j)$ fluid.
The inequality (\ref{eq:flFesibility}) enforces the physical constraint that
any given station can spend at most 100\% of its time processing
fluid.

Equations
(\ref{eq:flStage1}) and (\ref{eq:flStagej}) are simply flow
conservation equations: all class $(i,j-1)$ fluid becomes
class $(i,j)$ fluid after processing, for all $j\leq J_i$, and
class $(i,J_i)$ fluid leaves the network after processing. The
last constraint (\ref{eq:flNonnegativity}) simply says that all
the variables involved are non-negative real numbers.
Note that only the expectations
$1/\lambda_i= \mathbb{E}[\mX^i_1]$ and
$1/\mu_{i,j}=\mathbb{E}[\mS^{i,j}_1]$ of interarrival and service
times appear in the fluid model. The higher order moments
of the network primitives are not reflected in the model.

For each station $\sigma$ we let \be{eq:flQsigma}
\bQ_{\sigma}(t)=\sum_{(i,j)\in\sigma}\bQ_{i,j}(t), \ee that is,
$\bQ_{\sigma}(t)$ is the total fluid level in station $\sigma$ at
time $t$. Also let \be{eq:flTsigma}
\bT_{\sigma}(t)=\sum_{(i,j)\in\sigma}\bT_{i,j}(t). \ee So,
$\bT_{\sigma}(t)$ is the total amount of time station $\sigma$
spent processing  fluid during the time interval $[0,t]$.
Equivalently, $\bI_{\sigma}(t)\equiv t- \bT_{\sigma}(t)$ represents the
cumulative amount of idling experienced by  station $\sigma$
during the time interval $[0,t]$.

>From  inequality (\ref{eq:flFesibility}) it follows that the
function $\bT_{i,j}(t)$ is Lipschitz  continuous. Using
equations (\ref{eq:flStage1})--(\ref{eq:flNonnegativity}) it
can be checked that all of
$\bQ_{i,j}(t),\bA_i(t),\bD_{i,j}(t),\bT_{i,j}(t)$ are also
Lipschitz continuous. Any solution
$(\bQ_{i,j}(t),\bA_i(t),\bD_{i,j}(t),\bT_{i,j}(t))$ of the system
of equations and inequalities
(\ref{eq:flStage1})--(\ref{eq:flNonnegativity}) is defined to be a
fluid solution.
For simplicity, henceforth we use $(\bQ(t),\bT(t))$ to denote a fluid solution, where
$\bQ(t)$ and $\bT(t)$ stand respectively for vectors
$(\bQ_{i,j}(t))$ and $(\bT_{i,j}(t))$. A fluid solution
$(\bQ(t),\bT(t))$ is defined to be non-idling if for every station
$\sigma$, $\bI_{\sigma}(t)$ increases only at times $t$ when
$\bQ_{\sigma}(t)=0$. Formally, the fluid solution is non-idling if
for every station $\sigma$
\be{eq:flNon-idling}
\int_0^{\infty}\bQ_{\sigma}(t)d\bI_{\sigma}(t)=0.
\ee
The integral
is well-defined because $\bI_{\sigma}(t)$ is a Lipschitz continuous
function and, as a result, is almost everywhere differentiable in
$\Re_+$ with respect to the Lebesgue measure on $\Re_+$.

\begin{Defi}\label{def:FluidModel}
The set of non-idling feasible solutions to the system of
equations~(\ref{eq:flStage1})-(\ref{eq:flNonnegativity}),(\ref{eq:flNon-idling})
is defined to be the non-idling \emph{fluid model}.
\end{Defi}

When a queueing network operates under a specific scheduling policy, for example,
under a fixed buffer priority policy, additional constraints can be
added to the fluid equations in order to reflect the policy. In this
paper we are only considering the case of all the non-idling policies,
and thus the non-idling fluid model defined is the one of interest.
For the remainder of the paper we drop the modifier ``non-idling''
and simply refer to the ``fluid model.''

The following lemmas are used later in the paper. The proofs of
both lemmas are straightforward and thus omitted (note that
Lemma \ref{lemma:scaling} appeared as Property 2 in Chen~\cite{che95}).

\begin{lemma}\label{lemma:linearization}
Suppose $(\bQ(t),\bT(t))$ is a fluid solution defined over a time interval $[0,\theta]$.
Then $\bQ'(t) \equiv \bQ(0)+{t\over \theta}(\bQ(\theta)-\bQ(0)), \bT'(t) \equiv {t\over \theta}\bT(\theta)$
defined over $[0,\theta]$
is also a fluid solution. Moreover, suppose the solution $(\bQ(t),\bT(t))$ is non-idling and
for every station $\sigma$ either
$\bQ_{\sigma}(t)>0$ for all $t\in [0,\theta]$ or $\bQ_{\sigma}(0)=\bQ_{\sigma}(\theta)=0$.
Then the solution $(\bQ'(t),\bT'(t))$ is also non-idling.
\end{lemma}

\begin{lemma} \label{lemma:scaling}
Suppose $(\bQ(t),\bT(t))$ is a non-idling fluid solution defined over a time interval $[0,\theta]$.
Then for any $\beta > 0$, $\bQ'(t) \equiv \beta \bQ(\beta^{-1}t)$,
$ \bT'(t) \equiv \beta \bT(\beta^{-1}t)$ is a non-idling fluid solution defined
over the interval $[0, \beta \theta]$.
\end{lemma}

In the proofs in later sections, we need to define certain
types of fluid models
with a finite decomposition property. We define this notion
below.

\begin{Defi}\label{Defi:decomposition}
A fluid model  is defined to satisfy the \emph{Finite
Decomposition Property (FDP)} if there exist values $\nu, B>0$, with the
following property. For every  non-idling
fluid solution $(\bQ(t),\bT(t))$
defined over an interval $[0,\theta]$
such that $\bQ(t) \not= 0$ on this interval,
there exist a non-idling
fluid solution $(\tilde Q(t),\tilde T(t))$ also defined over
$[0,\theta]$ and a sequence of times instances
$0=t_0<t_1<t_2<\cdots<t_M=\theta$ such that
\begin{enumerate}
\item $M\leq \nu\theta\sup_{0\leq t\leq \theta}{1\over ||\bQ(t)||} + B$ and
$\inf_{0\leq t\leq \theta} ||\tilde Q(t)||\geq \inf_{0\leq t\leq \theta} ||\bQ(t)||.$
\item $\tilde Q(t_m)=\bQ(t_m)$ for all $m=0,1,\ldots,M$.
\item For each interval $(t_r,t_{r+1}), 0\leq r\leq M-1$ and each station $\sigma$ either
$\tilde Q_{\sigma}(t)>0$ for all $t\in(t_r,t_{r+1})$ or $\tilde
Q_{\sigma}(t)=0$ for all $t\in(t_r,t_{r+1}).$
\end{enumerate}
\end{Defi}

The next proposition shows that the FDP requirement is not
restrictive for fluid models arising from two station networks.

\begin{prop}\label{prop:FiniteDecomposition}
Fluid networks with two stations ($J=2$) satisfy  FDP.
\end{prop}

Although we only consider multitype fluid networks in this paper,
the proposition actually holds for any two station fluid network, for example networks with proportional routing.
This general form of Proposition  \ref{prop:FiniteDecomposition} is proved in
Subsection \ref{subsection:FDJ=2}.
At this point we do not know whether FDP holds for
general networks (i.e. with $J >2$).

\subsubsection{Global stability and global weak stability}\label{subsubsection:StabilityFluid}
Just as for stochastic queueing networks, we can define
stability and global stability for fluid networks.

\begin{Defi}\label{Defi:flStable}
A fluid solution $(\bQ(t),\bT(t))$ is defined to be stable if
there exists a $\tau<\infty$ such that $\bQ(t)=0$ for all
$t\geq\tau$. A fluid model is defined to be globally stable if
there exists a $\tau<\infty$ such that  every non-idling fluid
solution $(\bQ(t),\bT(t))$ satisfying $||\bQ(0)||=1$ also
satisfies $\bQ(t)=0$ for all $t\geq\tau$.
\end{Defi}

\remarks

1. The condition $||\bQ(0)||=1$ in the definition above is a
necessary scaling condition. One cannot have a uniform emptying
time $\tau$ without a bound on the initial state.

2. The definition of  global stability is somewhat different from
the perhaps more natural: ``network is defined to be globally stable
if it is stable for all non-idling policies.'' While it is possible
that both definitions are equivalent and it is known to hold
in many cases, it has not yet been established in general.
Definition \ref{Defi:flStable} is used more often because
it simplifies certain technical considerations.

Below, we define a stability notion for fluid networks which is
the analogue of the rate stability definition for stochastic
networks.

\begin{Defi}\label{Defi:flWeakStable}
A fluid model is defined to be globally weakly stable if for any
non-idling fluid solution $(\bQ(t),\bT(t))$, $\bQ(0)=0$ implies
$\bQ(t)=0$ for all $t \ge 0$.
\end{Defi}

In words, a fluid model is weakly stable if one cannot construct a
non-zero fluid solution which starts from zero. We did not
introduce the notion of a weakly stable fluid solution, since this
would just mean introducing a trivial $\bQ(t)=0$ solution (also it is
easy to check that $\bQ(t)=0$ for all $t$ implies
$T_{i,j}(t)={\lambda_i\over \mu_i}t$ for all $(i,j)$).

\subsection{The connections between stochastic and fluid queueing networks}\label{subsection:stoch fluid}

The most immediate connection between a stochastic network and the
corresponding fluid queueing network is provided by the results of
Dai \cite{dai95b} and Stolyar \cite{sto95}. Roughly speaking,
they show that for a broad class of scheduling policies, if
a stochastic network is operating under a policy
${\cal U}$, each weak limit $\bQ(t)=\lim_n {\mQ(nt)\over t}$
and
$\bT(t)=\lim_n {\mT(nt)\over t}$ of the stochastic queue length
process $\mQ(t)$ and cumulative work process $\mT(t)$, with a
sequence of initial states $\mQ(0)=\lfloor n\gamma\rfloor$, where
$\gamma$ is a fixed positive constant, is a deterministic
continuous function $(\bQ(t),\bT(t))$ which is a fluid solution of
the corresponding fluid model. If the policy ${\cal U}$ is
non-idling, then each obtained fluid solution is also non-idling.
Thus, the  queue length process after an appropriate rescaling
using certain scaled initial states, converges to a fluid solution.

This rescaling  process provides the basic tool for
connecting the stability of stochastic and fluid networks. In
fact, this connection was the primary motivation for introducing
fluid model techniques \cite{rybsto92}. The following theorem
establishes a fundamental relationship between the stability of
the stochastic and fluid models.

\begin{theorem}\label{theorem:DaiStolyar}(Dai \cite{dai95b}, Stolyar \cite{sto95})
Consider a multitype queueing network. If the corresponding
fluid model is globally stable then the stochastic network is globally stable.
\end{theorem}

The theorem actually holds for a broader class of networks and
also for networks operating under specific scheduling policies.
If one is given a particular scheduling policy ${\cal U}$, one can
sometimes identify additional constraints that the fluid limits
$\lim_n \mQ(nt)/n$ must satisfy.

More relevant to the topic of the present paper is the following
related result.

\begin{theorem}\label{theorem:Chen}(Chen \cite{che95})
Consider a multitype queueing network. If the corresponding
fluid model is globally weakly stable then the
stochastic network is globally rate stable.
\end{theorem}

Our understanding of global stability and global weak stability is
fairly complete for fluid models corresponding to queueing
networks with two stations ($J=2$), thanks to the results of
Bertsimas et al. \cite{bgt96} and Dai and Vande Vate \cite{daivan97}.
Both of these papers obtain necessary and sufficient conditions
for global stability of fluid networks for the case $J=2$.
Moreover a certain parameter $\rho^*$ is introduced in
\cite{daivan97}. This parameter is called  the maximum virtual traffic
intensity.
It  is shown that the fluid
model is globally stable iff  $\rho^*<1$ and is globally weakly
stable iff $\rho^*\leq 1$. The condition $\rho^*\leq 1$ then implies
rate stability of the underlying stochastic network by Theorem
\ref{theorem:Chen}. One of the main results of our paper is to
establish a converse: $\rho^*>1$ implies the stochastic network is
not globally rate stable. In particular, $\rho^*\leq 1$ is the
tight global rate stability condition for multitype networks
with two stations.

\section{Main Results}\label{section:main}

In this section, we provide the main results and corollaries of
this paper. All proofs, along with the needed lemmas, are provided
in Section~\ref{sec:proofs}. Our first result concerns the
structural properties of non-idling fluid solutions. The result is
introduced primarily because it is needed to prove the main result
of the paper, but we believe that it is interesting in its own
right and thus state the result in this section.

\begin{theorem}\label{theorem:FluidLinear}
Suppose the fluid model of a multitype queueing network is not weakly stable. Then there
exists a positive constant $\gamma>0$  such that for any initial
state $q\in\Re_+^d$, there exists a non-idling fluid solution
$(\bQ(t),\bT(t))$ satisfying $\bQ(0)=q$ and $||\bQ(t)||\geq \gamma t$ for
all $t\geq 0$. Namely, the solution is linearly divergent.
Moreover, this solution satisfies
\be{eq:lowerFluid} \inf_{t\geq 0}||\bQ(t)||\geq {||q||\over
2}\min\left({\gamma\over C},1\right), \ee
where $C$ is defined by (\ref{eq:C}).
\end{theorem}

Intuitively, the notion of a fluid model not being weakly stable
seems weaker than linear divergence. In particular, a fluid model is
not weakly stable if there exists a solution which ``pops up from zero''
at some point, after starting in the zero state.
Theorem~\ref{theorem:FluidLinear} shows that if one solution pops up,
then a different solution can be constructed which
diverges to infinity linearly, i.e.\ we construct a stronger fluid solution
(in the sense of instability) from a seemingly weaker solution.
This stronger fluid solution can then be used to infer the
instability of a class of associated stochastic networks. Finally,
we note that the divergent solution can be constructed from any
initial state $\bQ(0)=q\in \Re_+^d$, including the zero state.

We are now prepared to state the main result of the paper, which
connects the instability of fluid models and stochastic networks.

\begin{theorem}\label{theorem:Main}
Consider a  multitype stochastic network
satisfying Assumptions A, B, and C. Suppose the associated fluid
model is not globally weakly stable, and satisfies FDP.
Then, for any  initial state $(q,z_1,z_2)\in \mathbb{Z}_+^d\times \Re^{I+d}$, there exists a non-idling scheduling
policy for which the resulting queue level process satisfies
$\liminf_{t\rightarrow\infty}{||\mb{Q}(t)||\over t}>0$ with
probability one. In particular, the stochastic process
associated with the queueing network is unstable, under some non-idling policy.
\end{theorem}
The rate of divergence to infinity implied by the theorem above will be explicit. We will show
that constructed policy results in
\be{eq:rate0}
\liminf_{t\rightarrow\infty}{||\mb{Q}(t)||\over t}\geq {\max\left({\gamma\over C},1\right)\over 8\max(1,{3\over \gamma})}
\ee
with probability one.
The import of Theorem~\ref{theorem:Main} is more apparent from the
corollaries provided below.

\begin{coro} \label{coro:1}
Consider a multitype stochastic network with $J=2$. If the associated fluid
model is not globally weakly stable, then the queueing network is
unstable in the sense that
$\liminf_{t\rightarrow\infty}{||\mb{Q}(t)||\over t}>0$ with
probability one from each initial state under some non-idling
scheduling policy.
\end{coro}

Corollary~\ref{coro:1} follows from Theorem~\ref{theorem:Main} and
Proposition \ref{prop:FiniteDecomposition}, which states that
FDP holds for fluid networks with two stations. Recall that one motivation for
our work is the stability Theorems~\ref{theorem:DaiStolyar} and \ref{theorem:Chen}.
Thus Corollary~\ref{coro:1} provides a complete converse of
Theorem~\ref{theorem:Chen} for two station multitype networks.
A missing piece in the theory for general $J$ is to determine if all
fluid models satisfy FDP. If such a result holds, then
Theorem~\ref{theorem:Main} would imply a converse for networks with
an arbitrary number of stations.

We note also that Theorem~\ref{theorem:Chen} is valid when we consider
fluid and queueing networks under specific scheduling
policies. However, for networks operating under specific policies
(rather than a class of policies) a general converse to the
theorems of Chen and Dai is not possible as demonstrated in
\cite{dhv00}.

Dai and VandeVate~\cite{daivan97} derived explicit necessary and
sufficient conditions for global weak stability of fluid models of multitype
networks in terms of a
certain parameter $\rho^*$ related to the so-called
virtual traffic intensity and push start conditions.
They prove that such fluid networks are weakly stable iff $\rho^*\leq 1$.
Considering Theorem~\ref{theorem:Chen} along
with Corollary~\ref{coro:1}, those results now yield complete
necessary and sufficient conditions for rate stability of two
station stochastic mutlitype networks.

\begin{coro} \label{coro:2}
A stochastic two station multitype network  is globally rate stable if and only
if $\rho^*\leq 1$.
\end{coro}

\section{Proofs of Main Results} \label{sec:proofs}

In this section we provide all of the proofs of our main results.
The first proof, of Theorem \ref{theorem:FluidLinear}, shows that
if the fluid is not globally weakly stable, there exists a linearly
divergent fluid solution.

\subsection{Linearly divergent fluid solutions}\label{subsection:fluid}

\proof[Proof of Theorem \ref{theorem:FluidLinear}]
We assume in the theorem that a given fluid model is not globally weakly stable. Hence, there
exists a non-idling solution which satisfies $\bQ(0)=0$ and $\bQ(t_0) \not= 0$
for some $t_0 > 0$. First note that, without the loss of generality, we may assume
that $\bQ(t)\neq 0$ for all $0<t\leq t_0$. Otherwise,
we can find $\hat{t} = \sup \{ 0 \le t < t_0 : \bQ(t) = 0 \}$
and consider the fluid solution on $[\hat{t}, t_0]$ only.
Note that $\hat{t} < t_0$ by the continuity of $\bQ(t)$ and
the fact that $\bQ(t_0) \not= 0$. Next, using Lemma \ref{lemma:scaling}
with some $\beta > 0$ we can obtain a new solution defined
on $[0, \beta t_0]$ with $\bQ'(0)=0$ and $||\bQ'(\beta t_0)||=\beta ||\bQ(t_0)||$.
If we set $\beta = t_0^{-1}$ then we have a solution defined
on $[0,1]$ with $\bQ'(0)=0$ and $||\bQ'(1)||= ||\bQ(t_0)|| / t_0$.
Hence, again without loss of generality, we set $t_0=1$, i.e.\ we assume we
are given
a non-idling solution with $\bQ(0)=0$ and $\bQ(1) \not= 0$.

We now build a new fluid solution by constructing it iteratively
over the intervals $[0,1)$, $[1,2)$, $[2,4)$, \ldots,
$[2^n,2^{n+1}), \ldots.$ We denote the solution that
is constructed in this manner by
$(\bQ^{o}(t),\bT^{o}(t))$.
For the initial interval $[0,1)$,
consider our initial fluid solution $\bQ(t)$ satisfying $\bQ(0)=0,\bQ(1)\neq 0$.
We first modify the solution by setting $\bQ(0) = q$, where
$q \in \Re_+^d$. Next for
every $t\leq 1$ and every class $(i,j)$,
on the interval $[0,t]$ we spend exactly $\bT_{i,j}(t)$
time units processing class $(i,j)$ flow,
\emph{plus} whatever necessary additional amount is required to
make the solution non-idling. In other words, we can think of the
flow ``created'' and ``processed'' by the non-weakly stable solution
$\bQ(t)$ as high priority flow, and the remaining flow as low priority flow.
Note that the allocation of the additional processing effort required
is not necessarily uniquely determined by the original allocation $\bT(t)$.
In any case, the resulting solution
satisfies $\bQ^{o}(0)=q$ and $\bQ_{i,j}^{o}(t)\geq\bQ_{i,j}(t)$ for
all classes $(i,j)$ and $t\leq 1$. In particular, $||\bQ^{o}(1)||\geq||\bQ(1)||>0$.

Assume now the solution has been constructed over the time horizon $[0,2^n]$ for $n\geq 0$.
We now extend it over $[2^n,2^{n+1}]$.
The idea of the construction is similar to the first interval, except that
we ``stretch'' the original solution $\bQ(t)$ by a factor of $2^n$ and
then use this solution to extend our current solution by defining it on
$[2^n, 2^{n+1}]$. That is, consider
the scaled solution $(\beta\bQ(\beta^{-1} t),\beta\bT(\beta^{-1} t))$ with $\beta=2^n$.
This solution is defined
over $t\in [0,2^n)$. Next, for each $t\in [2^n,2^{n+1}]$ let $\bT^{o}(t)$ be defined by
$\bT^{o}(t)-\bT^{o}(2^n)=2^n\bT(2^{-n}(t-2^n))$ plus any extra processing
effort required to make the solution non-idling.

It can be easily checked  that the resulting solution $\bQ^{o}(t)$ satisfies $\bQ^{o}_{i,j}(t)\geq 2^{n}\bQ_{i,j}
(2^{-n}(t-2^n ))$
for all $t\in [2^n,2^{n+1}]$ and all $i,j$
which implies $||\bQ^{o}(t)||\geq 2^{n}||\bQ(2^{-n}(t-2^n))||$.
In particular, $||\bQ^{o}(2^{n+1})||\geq 2^{n }||\bQ(1)||$.

We have constructed a non-idling fluid solution $\bQ^{o}(t)$ which diverges
to infinity at time instances $t_n=2^n, n=0,1,\ldots.$
To complete the proof of the theorem, we show that for some constant
$\gamma_0>0$, $||\bQ^{o}(t)||\geq \gamma_0 2^n$
for all $t \in [2^n,2^{n+1}]$.
First let us show that this implies the theorem. For any $t >0$ find
the largest integer $n$ such that $2^n\leq t$, i.e.\
let $n=\lfloor\log_2 t\rfloor$. We have $||\bQ^{o}(t)||\geq \gamma_0 2^n\geq
\gamma_0 2^{\log_2 t-1}=\gamma_0 t/2.$ Setting $\gamma=\gamma_0/2$, we obtain the result.

To show the existence of $\gamma_0$, note that for any $t_1<t_2$ and any feasible
fluid solution $\bQ(\cdot)$ we have:
\begin{eqnarray}
||\bQ(t_2)|| & \geq & ||\bQ(t_1)||-\sum_{1\leq i\leq I}\mu_{i,J_i}(t_2-t_1) \notag\\
& \geq & ||\bQ(t_1)||-C(t_2-t_1). \label{eq:-Ct}
\end{eqnarray}
This implies that for all $ t \in [ 2^n,2^n+2^{n-1}||\bQ(1)||/(2C)]$, $\bQ^{o}(t)$
satisfies:
\begin{eqnarray*}
||\bQ^{o}(t)|| & \geq  & ||\bQ^{o}(2^n)||-C(t-2^n) \\
& \geq & 2^{n-1}||\bQ(1)||-C(t-2^n) \\
& \geq & 2^{n-2}||\bQ(1)||.
\end{eqnarray*}
 If $2^n+2^{n-1}||\bQ(1)||/(2C)\geq 2^{n+1}$,
then we simply set $\gamma_0=(1/4)||\bQ(1)||$. Otherwise, let
$$\gamma_1=\min \left \{||\bQ(t)||:{||\bQ(1)||\over 4C}\leq t\leq 1
\right \}.$$
This minimum exists since $\bQ(t)$ is continuous and
it is positive since $|| \bQ(t)|| > 0$ for all $0<t\leq 1$. Then, for all
$2^n+2^{n-1}||\bQ(1)||/(2C)\leq t\leq2^{n+1}$, we have $||\bQ^{o}(t)||\geq 2^n||\bQ(2^{-n}(t-2^n))||
\geq 2^n\gamma_1$. We take $\gamma_0=\min\{(1/4)||\bQ(1)||,\gamma_1\}$ and we
have proven the first inequality in the theorem statement.

%%DG
The last part of the proposition follows almost immediately. Using (\ref{eq:-Ct})
with $t_1 = 0$ and $t_2 = t$ we have
$||\bQ(t)||\geq ||q||-Ct\geq ||q||/2$ for $t\leq ||q||/(2C)$.
On the other hand, by construction, $||\bQ(t)||\geq \gamma t\geq \gamma ||q||/(2C)$, whenever
$t\geq ||q||/(2C)$. This completes the proof of the theorem. \qed
%%EndDG

Theorem~\ref{theorem:FluidLinear} will  be used for proving our main result, Theorem \ref{theorem:Main}.
%we will follow the same line of arguing when proving Theorem~\ref{theorem:Main}.
Specifically, we will construct a non-idling scheduling policy for the
discrete network which, with high
probability, results in a trajectory very close to the fluid trajectory built
in the proof of Theorem~\ref{theorem:FluidLinear}.
We will use the large deviations bounds (\ref{eq:LDBtime}) and (\ref{eq:LDBrate}) multiple
times to obtain bounds on the deviation between the fluid and stochastic
trajectories.

\subsection{FDP in fluid networks with two stations}\label{subsection:FDJ=2}

\proof[Proof of Proposition~\ref{prop:FiniteDecomposition}]
Consider
a network with two stations, $\sigma_1$ and $\sigma_2$, and suppose we have
a non-idling fluid solution $(\bQ(t), \bT(t))$ which is non-zero
over time
interval $[0,\theta]$. By continuity, $\inf_{0\leq t\leq
\theta}||\bQ(t)||>0$. The next result follows from
Proposition 1 in \cite{bgt96}.
There exists a nondecreasing sequence $t_i$
such that $\sup_i t_i=\theta$ and such that for all times less than $\theta$
the following hold:
\begin{itemize}
\item $\bQ_{\sigma_1}(t_{4m+1})>0,\bQ_{\sigma_2}(t_{4m+1})=0$ and
for  $ t\in [t_{4m+1},t_{4m+2}],~  \bQ_{\sigma_1}(t)>0$;
\item $\bQ_{\sigma_1}(t_{4m+2})>0,\bQ_{\sigma_2}(t_{4m+2})=0$ and
for  $ t\in (t_{4m+2},t_{4m+3}), ~
\bQ_{\sigma_1}(t),\bQ_{\sigma_2}(t)>0$;
\item $\bQ_{\sigma_2}(t_{4m+3})>0,\bQ_{\sigma_1}(t_{4m+3})=0$ and
for  $ t\in [t_{4m+3},t_{4m+4}],~\bQ_{\sigma_2}(t)>0  $;
\item $\bQ_{\sigma_2}(t_{4m+4})>0,\bQ_{\sigma_1}(t_{4m+4})=0$ and
for  $ t\in (t_{4m+4},t_{4m+5}),~\bQ_{\sigma_1}(t),\bQ_{\sigma_2}(t)>0.$
\end{itemize}

Moreover, one of $t_i, i=1,2,3,4$ is equal to zero. When $t_2,t_3$ or $t_4$
is zero,
$t_i$ with lower value of $i$ is not defined.

The characterization above essentially divides the trajectory of a fluid
solution into four different segments. On the segment of the trajectory
between $t_{4m+1}$ and $t_{4m+2}$, the trajectory is either on the boundary
of the state space (where $\bQ_{\sigma_2}(t) = 0$) or in the interior
of the state space. We next claim that such a segment can be ``linearized''
such that it remains a non-idling solution, yet $\bQ_{\sigma_2}(t) = 0$
for all $t \in [t_{4m+1},t_{4m+2}]$. In other words, the linearized solution
is on the boundary for the entire interval. To achieve the linearization
we define
\begin{eqnarray*}
\tilde{Q}(t) & = & \bQ(t_{4m+1}) + \frac{t-t_{4m+1}}{t_{4m+2}-t_{4m+1}}
[ Q(t_{4m+2}) - Q(t_{4m+1})] \quad \mbox{and} \\
\tilde{T}(t) & = & \bT(t_{4m+1}) + \frac{t-t_{4m+1}}{t_{4m+2}-t_{4m+1}}
[ T(t_{4m+2}) - T(t_{4m+1})],
\end{eqnarray*}
for all $t \in [t_{4m+1},t_{4m+2}]$.
Using Lemma~\ref{lemma:linearization} it follows that the new solution $(\tilde{Q}(t),\tilde{T}(t))$
is both feasible and non-idling, given that the original solution was also.
In a similar manner, we linearize the fluid solution $(\bQ(t),\bT(t))$
on all intervals of the form
$[t_{4m+3}, t_{4m+4}]$. Hence, in each interval the new solution remains on
one of the axes,
unless it is crossing the interior, from one axis to the other.

We now demonstrate that $(\tilde{Q}(t),\tilde{T}(t))$ has the
properties described in Definition~\ref{Defi:decomposition}.
First, we claim that for each $m$,
\be{eq:t4m}
t_{4m+3}-t_{4m+1}\geq \inf_{0\leq t\leq \theta}||\bQ(t)||/C, \qquad
t_{4m+3}-t_{4(m+1)+1}\geq \inf_{0\leq t\leq \theta}||\bQ(t)||/C.
\ee
Indeed, by construction
$\tilde{Q}_{\sigma_1}(t_{4m+1})>0,\tilde{Q}_{\sigma_2}(t_{4m+1})=0$
and $\tilde{Q}_{\sigma_2}(t_{4m+3})>0,\tilde{Q}_{\sigma_1}(t_{4m+3})=0$.
In particular
$$\tilde{Q}_{\sigma_1}(t_{4m+1})=||\tilde{Q}(t_{4m+1})||\geq
\inf_{0\leq t\leq \theta}||\tilde{Q}(t)|| \geq \inf_{0\leq t\leq
\theta}||\bQ(t)|| > 0.$$
Note that total rate at which fluid can depart from a
given station is bounded above by
$\sum_{i,j}\mu_{i,j}<C$. Thus, since
$\tilde{Q}_{\sigma_1}(t_{4m+3})=0$, we have $t_{4m+3}-t_{4m+1}\geq \inf_{0\leq
t\leq
\theta}||\bQ(t)||/C$. An analogous argument demonstrates
that $t_{4m+3}-t_{4(m+1)+1}\geq \inf_{0\leq t\leq
\theta}||\bQ(t)||/C$.
Since the interval lengths are bounded strictly away from zero,
the total number of
points $t_i$ in $[0,\theta]$ is at most $(2\theta C/\inf_{0\leq t\leq
\theta}||\bQ(t)||)+2$, where the $+2$ accounts
for the end points of $[0,\theta]$.
Setting $\nu=2C, B=2$ yields the first FDP property (1).
Properties (2) and (3) are automatically satisfied by our construction
of $(\tilde{Q}(t),\tilde{T}(t))$ above.
\qed

\subsection{Transient paths in the stochastic network}\label{section:stochastic}

Most of this subsection is devoted to the proof of Theorem~\ref{theorem:mainStoch}
of Section \ref{sec:mp},
which, as we will show, implies the main result of our paper, Theorem~\ref{theorem:Main}.
In the proof we repeatedly use probabilistic bounds of the form $c_1\exp(-c_2n)$,
where $c_1,c_2> 0$ are  constants which depend on the parameters of our queueing network
and $n$ is a scaling parameter which takes on
a large value. In various expressions,
$c_2$ is usually related to the constant $L$ appearing in the large deviations
bounds in (\ref{eq:LDBtime}) and (\ref{eq:LDBrate}) and the network parameters
$\lambda_i,\mu_{i,j},|I|,C,$ as well as parameter $\gamma$ introduced in Theorem \ref{theorem:FluidLinear}.
We will also be considering finite sums of the bounds
of the form $c_1'\exp(-c_1n)+c_2'\exp(-c_2 n)+...+c_m'\exp(-c_m n)$. In general
the $c_i,c_i'$ take on different values and $m$ is a constant, independent
of $n$. Such sums can be bounded above by $c'\exp(-c n)$ for
$c=\min_{1\leq k\leq m}c_k$ and $c'=\sum c_i'$.

In our proofs, the actual values of the constants are not important,
only the fact that they are independent of $n$. Therefore,
to simplify the exposition, we simply use the notation $O(\exp(-\Theta(n)))$ and we
write expressions like
$O(\exp(-\Theta(n)))+O(\exp(-\Theta(n)))=O(\exp(-\Theta(n)))$ where the standard notations $O(\cdot)$
and $\Theta(\cdot)$ hide the actual constants $c$ and $c'$.

\subsubsection{Proof Preliminaries and the Scheduling Policy $\cal{U}$}

In order to precisely state the next series of detailed results, we need
to define a non-idling policy $\cal{U}$. The definition of this policy
involves a number of preliminary observations and definitions.

First, let \be{eq:bound_theta}
\theta=\max\left(1,{3\over \gamma}\right).
\ee
The parameter $\theta$ depends only on parameters of the model since $\gamma$ depends only
on the parameters of the model.

Consider any initial state $(q,z_1,z_2)\in \mathbb{Z}_+^d\times \Re_+^{I+d}$.
Let $n=||q||$. By Theorem~\ref{theorem:FluidLinear} there
exists a non-idling fluid solution $(\bQ(t),\bT(t))$ which satisfies $\bQ(0)=q$ and
$||\bQ(t)||\geq\gamma t$ for all $t\geq 0$.
Since FDP is assumed then by Proposition~\ref{prop:FiniteDecomposition}
the solution $(\bQ(t),\bT(t))$ can be modified to a solution which satisfies properties
described in Definition~\ref{Defi:decomposition}.
Let
\be{eq:theta=2gamma}
\theta_0=\theta ||q||=||q||\max\left(1,{3\over \gamma}\right),
\ee
in which case we have
\be{eq:bQ>3}
||\bQ(\theta_0)||\geq 3||q||.
\ee
By Theorem~\ref{theorem:FluidLinear}, the fluid solution
is also such that
\be{eq:infQ}
\inf_t||\bQ(t)||\geq {||q||\over 2}\min\left({\gamma\over C},1 \right).
\ee
Since FDP is satisfied, there exists
another solution $(\tilde Q(t),\tilde T(t))$ and a sequence $0=s_0<s_1<\cdots<s_M=\theta_0$, such that
$\inf_{0\leq t\leq \theta_0}||\tilde{Q}(t)||\geq\inf_{0\leq t\leq \theta_0}||\bQ(t)||$ and for every
interval $[s_r,s_{r+1}]$ and for each station $\sigma$ either $\tilde Q_{\sigma}(t)$ is zero within
$(s_r,s_{r+1})$, or it is strictly positive within $(s_r,s_{r+1})$.
For simplicity we assume that $(\bQ(t),\bT(t))$ is this modified solution.
In such a modified solution we also note that
\begin{eqnarray}
M\leq \nu\theta_0\sup_{0\leq t\leq \theta_0}{1\over ||\tilde{Q}(t)||}+B & \leq &
\nu\theta_0 {2\over ||q||}\max\left({C\over \gamma},1\right)+B \notag \\
& \leq & 2\nu\max\left(1,{3\over \gamma}\right)\max\left({C\over \gamma},1\right)+B \label{eqM<}
\end{eqnarray}
where we used (\ref{eq:theta=2gamma}) and (\ref{eq:infQ}).
In particular, we obtain a bound
on $M$ which depends  only on the parameters of the model (and is independent of $||q||$),
since $\nu,\gamma$ and $C$ depend only on the parameters of the model. Note, on the other hand, that
the partition $s_r, r=0,1,\ldots,M$ does depend on $q$.
Recalling the notation $||q||=n$ we rewrite (\ref{eq:bQ>3}) and (\ref{eq:infQ}) as
\be{eq:bQn>3}
||\bQ(n\theta)||\geq 3||q||=3n,
\ee
and
\be{eq:bQn1>3}
\inf_{0\leq t\leq \theta n}||\bQ(t)||\geq {n\over 2}\min\left({\gamma\over C},1\right).
\ee
Our next goal is to describe a non-idling scheduling policy ${\cal U}={\cal U}(\delta)$ implemented
over the time horizon $[0,\theta_0]=[0,\theta n]$. Recall that our starting state is $(q,z_1,z_2)$.
In particular, $\mQ(0)=\bQ(0)=q$. The policy ${\cal U}$ attempts to mimic
the fluid solution described above, over the same time interval.
We parameterize the policy with a constant $\delta >0$ which is any constant satisfying
\be{eq:delta}
\delta\leq {1\over 12C^{M+3}}\min\left({\gamma\over C},1\right).
\ee
Let $t_m=m\delta n$ for $m=0,1,\ldots,\lceil\theta/\delta\rceil$.
We describe the policy ${\cal U}$ on each time interval $I_m=[t_m,t_{m+1})$.
For each time interval $I_m$ each station $\sigma$ nominally allocates $\bT_{i,j}(t_m,t_{m+1})$
time units to serving class $(i,j)$, for every class $(i,j)\in\sigma$.
To be precise, we first order all the classes at a station in a fixed,
but arbitrary manner. During the interval a class $(i,j)$ is chosen
for service, and we work on jobs from that class for
$\bT_{i,j}(t_m,t_{m+1})$ time units or until we exhaust the jobs
from class $(i,j)$.
Note that we cannot reach the end of the interval $I_m$, by the
feasibility of $\bT(\cdot)$ over this interval. When we are done
processing jobs of type $(i,j)$ the next class in the chosen order
of service is picked for processing. Note that we assume a preemptive
resume mechanism when switching between classes.
%% If a job is in service when the time allocated to a particular class
%% has been achieved, then the job is ejected and placed at the head
%% of the line for processing at a later time. When this class is again allocated time,
%% we assume that the remaining processing time for the ejected job
%% is the same as it was at the moment it was ejected, i.e.\ our policy is assumed
%% to be \emph{preemptive resume}.
%% It will be seen below, since the lengths of all intervals are order $\Theta(n)$, that the particular
%% preemption mechanism is irrelevant for the analysis and the argument goes through whether the policy
%% is preemptive or not.
If after going through all the classes, the time spent is strictly less
than $t_{m+1}-t_m$ and there are still jobs at the station, the station works on any available jobs.
If no jobs are available, the station idles.
Once the next time instance $t_{m+1}$ occurs, the policy is ``reset,''
in terms of the time allocations.

In other words, according to our scheduling policy, on each interval $I_m$ each station tries to spend
exactly the same amount of time on jobs in each class $(i,j)$
as the fluid solution $(\bQ(t), \bT(t))$ does, while maintaining the non-idling requirement.
Our main goal is to show that in general the resulting
stochastic process stays fairly close to the fluid trajectory $(\bQ(t),\bT(t))$, when
the stochastic network operates under the discipline ${\cal U}$.

>From the fluid equation $(\ref{eq:flFesibility})$ we have
$\sum_{(i,j)\in\sigma}\bT_{i,j}(t_m,t_{m+1})\leq t_{m+1}-t_m$, for
each $m$. As a result, any policy ${\cal U}$ is feasible. From the
description above, it is certainly non-idling. We now analyze the
dynamics of our network when policy ${\cal U}$ is implemented. For
convenience we introduce $s_{-1}\equiv s_0=0$.

\begin{lemma}\label{lemma:QQt_m}
Under the policy ${\cal U}$ (in fact under any scheduling policy), for every $m=0,1,\ldots,\lceil{\theta\over \delta}\rceil.$
\be{eq:Cdeltanbar}
\sup_{t_m\leq t\leq t_{m+1}}||\bQ(t)-\bQ(t_m)||\leq C\delta n.
\ee
and
\be{eq:mbQ2}
\mathbb{P}\left\{\sup_{t_m\leq t\leq t_{m+1}}||\mb{Q}(t)-\mb{Q}(t_m)|| > C\delta n\right\} \leq O(\exp(-\Theta(n))).
\ee
\end{lemma}

\proof
Applying (\ref{eq:flStage1}), (\ref{eq:flStagej}) and (\ref{eq:flFesibility}) we have
\[
||\bQ(t)-\bQ(t_m)||\leq (\sum_i\lambda_i+\sum_{i}\mu_{i,J_i})(t-t_m)<C(t_{m+1}-t_m)=C\delta n,
\]

which proves (\ref{eq:Cdeltanbar}). We now prove (\ref{eq:mbQ2}).
 By Assumption B (specifically bound (\ref{eq:LDBrate})), for every  $i$
and every $t \in [t_m, t_{m+1}]$
\begin{eqnarray}\label{eq:Ashort}
\mathbb{P}\left\{|\mb{A}_{i}(t)-\mb{A}_{i}(t_m)| > 2\lambda_i\delta n\right\} &\leq &
\mathbb{P}\left\{|\mb{A}_{i}(t_{m+1})-\mb{A}_{i}(t_m)| > 2\lambda_i\delta n\right\} \\
&\leq & O(\exp(-\Theta(n))),
\end{eqnarray}
since $t_{m+1}-t_m=\delta n$. Similarly, for all $i$ and $j$ and
$t \in [t_m, t_{m+1}]$
\be{eq:Dshort}
\mathbb{P}\left\{|\mb{D}_{i,j}(t)-\mb{D}_{i,j}(t_m)| > 2\mu_{i,j}\delta n\right\} \leq O(\exp(-\Theta(n))).
\ee
Applying (\ref{eq:stQueueStage1}) and (\ref{eq:stQueueStagej}) we obtain
\[
\mathbb{P}\left\{\sup_{t_m\leq t\leq t_{m+1}}|\mb{Q}_{i,j}(t)-\mb{Q}_{i,j}(t_m)| > 2(\lambda_i+\mu_{i,j-1}+\mu_{i,j})\delta n\right\} \leq O(\exp(-\Theta(n))).
\]
By summing these probabilities over all $(i,j)$, we obtain
\[
\mathbb{P}\left\{\sup_{t_m\leq t\leq t_{m+1}}||\mb{Q}(t)-\mb{Q}(t_m)|| > (2\sum_{i}\lambda_iJ_i+4\sum_{i,j}\mu_{i,j})\delta n\right\}
\leq \sum_{i,j}O(\exp(-\Theta(n)))=O(\exp(-\Theta(n))),
\]
implying
\[
\mathbb{P}\left\{\sup_{t_m\leq t\leq t_{m+1}}||\mb{Q}(t)-\mb{Q}(t_m)|| > C\delta n\right\} \leq O(\exp(-\Theta(n))),
\]
which is (\ref{eq:mbQ2}).
\qed

A large part of the remainder of the paper is devoted to proving Proposition
\ref{prop:StochasticClose} below. The proof is quite
lengthy and we split the argument into several subsections.

\begin{prop}\label{prop:StochasticClose}
Under the policy ${\cal U}={\cal U}(\delta)$ for every $r=-1,0,1,\ldots,M-1$, every $t_m\in [s_r,s_{r+1}]$ and every class $(i,j)$
\be{eq:st2}
\mathbb{P}\{|\mb{Q}_{i,j}(t_m)-\bQ_{i,j}(t_m)|\leq \delta C^{r+3} n\}\geq 1-O(\exp(-\Theta(n))).
\ee
\end{prop}

The proof is done by using various induction steps. The ``outer''
induction is on $r$, which indexes the trajectory decomposition
points $s_r$. The ``inner'' induction is done on the stages $j$ of the classes $(i,j)$ classes in
the network, and is outlined in various lemmas below.

We  start the outer induction with $r=-1$. Then for $t_m\in
[s_{-1},s_0]=\{0\}$ we simply have $t_m=0$ and the bound in
(\ref{eq:st2}) holds trivially for all classes $(i,j)$ since
$\mb{Q}(0)=\bar Q(0)=q,$ with probability one. Next we suppose the
bounds in (\ref{eq:st2}) hold for $-1,0,1, \ldots, r-1$. We then
show that the bounds hold for $r$. The necessary bounds will be
established by a sequence of lemmas. Our first lemma simply says that
assuming the bounds (\ref{eq:st2}) hold for all $r'\leq r-1$ and
$t_m\in [s_{r'},s_{r'+1}]$, a similar bound holds at the end point $s_r$.

\begin{lemma}\label{lemma:sr}
If the bound (\ref{eq:st2}) holds for all $r'\leq r-1$ then for every $i,j$
\[
\pr\left\{|\mb{Q}_{i,j}(s_{r})-\bQ_{i,j}(s_{r})|>\delta C^{r+2} n+2\delta C n \right\}\leq O(\exp(-\Theta(n))).
\]
\end{lemma}

\proof Find the largest $t_{m'}\leq s_r$. Then $t_{m'+1}=t_{m'}+\delta n > s_r\geq t_{m'}$.
By Lemma \ref{lemma:QQt_m} we have $|\bQ_{i,j}(s_{r})-\bQ_{i,j}(t_{m'})|\leq C \delta  n$ and
\[
\pr\left\{|\mb{Q}_{i,j}(s_{r})-\mb{Q}_{i,j}(t_{m'})|\geq C\delta n\right\}\leq O(\exp(-\Theta(n))).
\]
Since $t_{m'}\in [s_{r-1},s_r]$ then by the assumption of our induction in $r$,
\[
\pr\left\{|\mb{Q}_{i,j}(t_{m'})-\bQ_{i,j}(t_{m'})|\geq C^{r+2}\delta n\right\}\leq O(\exp(-\Theta(n))).
\]
Combining the last three inequalities, we obtain the result. \qed

In the next subsection we
obtain probabilistic lower bounds on the number of jobs processed
during the time interval $[s_r,t_m)$, for any
$t_m\in[s_r,s_{r+1}]$, under the scheduling policy ${\cal U(\delta)}$.

\subsubsection{Lower bounds on the departure process}\label{subsection:DepartureLower}

The next lemma shows that, with high probability, in the first stage in the route of each job
type the total number of jobs processed during the time interval  $[s_r,t_m]$
is not too far behind the amount of fluid processed during the same time interval in
the fluid solution. A subsequent lemma establishes a similar bound for stages two and higher.
Recall that we fixed $r$ and we assume by induction that (\ref{eq:st2}) holds for $r'\leq r-1$.

\begin{lemma}\label{lemma:DepartureStage1}
For every $i\leq I$ and every  $m$ such
that $s_r\leq t_m \leq s_{r+1}$,
\be{eq:st3}
\mathbb{P}\{\mb{D}_{i,1}(s_r,t_{m})\geq \bD_{i,1}(s_{r},t_{m})-2\delta C^{r+2} n\}\geq 1-O(\exp(-\Theta(n))),
\ee
and
\be{eq:st3time}
\mathbb{P}\{\mb{T}_{i,1}(s_{r},t_{m})\geq \bT_{i,1}(s_{r},t_{m})-3\mu_{\max}\delta C^{r+2} n\}\geq 1-O(\exp(-\Theta(n))).
\ee

\end{lemma}

\proof We start with proving bound (\ref{eq:st3}). Bound
(\ref{eq:st3time}) will be an easy corollary.

\textbf{Part I.} Fix a specific class $(i,1)$, and time
$t_{m_0},s_r\leq t_{m_0}\leq s_{r+1}$, and introduce the event
\be{eq:D}
{\cal D}(t_{m_0}) \equiv \left\{\mb{D}_{i,1}(s_r,t_{m_0}) < \bD_{i,1}(s_{r},t_{m_0})-2\delta C^{r+2} n \right\}.
\ee
Note then that (\ref{eq:st3}) is equivalent to having
$\pr\left\{D(t_{m_0})\right\}\leq O(\exp(-\Theta(n)))$ for every $i$ and $t_{m_0}\in [s_r,s_{r+1}]$.
Next, we introduce the events
\be{eq:A}
{\cal A} \equiv \left\{\forall \; t_m\in[s_r,s_{r+1}]: \mb{A}_{i}(s_r,t_m) \geq \lambda_i(t_m-s_r)-\delta C n \right\},
\ee
\be{eq:Q}
{\cal Q} \equiv \left\{\mb{Q}_{i,1}(s_{r})\geq \bQ_{i,1}(s_{r})-\delta C^{r+2} n-2\delta C n \right\}.
\ee
>From  Lemma \ref{lemma:sr} and the inductive assumption we have
\be{eq:Qprob}
\pr\left\{{\cal Q}\right\}\geq 1-O(\exp(-\Theta(n))).
\ee
Now fix any $t_m\in [s_r,s_{r+1}]$ and consider

$$\pr\{\mb{A}_{i}(s_r,t_m) \geq \lambda_i(t_m-s_r)-\delta C n \}=
\pr\left\{\mb{A}_{i}(s_r,t_m) \geq \lambda_i(t_m-s_r)-{\delta C n\over t_m-s_r}(t_m-s_r) \right\},$$
where without loss of generality we may assume $t_m>s_r$.
If $t_m-s_r\leq \delta n$, then the probability above is equal to one, since the right-hand side
of the inequality inside the probability is negative. Suppose now $t_m-s_r\geq \delta n$.
We have $${\delta C n\over t_m-s_r}\geq {\delta C n\over \theta n} = {\delta C\over \theta}.$$
Setting $\epsilon={\delta C\over \theta}$ and using the large deviations Assumption B with this
$\epsilon$ we obtain that
\be{eq:Alower}
\pr\{\mb{A}_{i}(s_r,t_m) \geq \lambda_i(t_m-s_r)-\delta C n \}\geq 1-O(\exp(-\Theta(t_m-s_r)))\geq 1-O(\exp(-\Theta(n))),
\ee
where $t_m-s_r\geq \delta n$ is used in the last inequality. The number of different $t_m$ in $[s_r,s_{r+1}]$ is at most $\theta n/(\delta n)=\theta/\delta$.
Summing over all such $t_m$ we conclude
\be{eq:AA}
\pr\{{\cal A}\}\geq 1-(\theta/\delta)O(\exp(-\Theta(n)))=1-O(\exp(-\Theta(n))).
\ee
Hence
\begin{eqnarray*}
\pr\{{\cal D}(t_{m_0})\} & = & \pr\{{\cal D}(t_{m_0})|{\cal A}\cap {\cal Q}\}\pr\{{\cal A}\cap {\cal Q}\}+
\pr\{{\cal D}(t_{m_0})|\overline{{\cal A}\cap {\cal Q}}\}\pr\{\overline{{\cal A}\cap {\cal Q}}\} \\
& \leq  &\pr\{D(t_{m_0})\cap A\cap Q\}+O(\exp(-\Theta(n))),
\end{eqnarray*}
where in the inequality we use $\pr\{\overline{{\cal A} \cap {\cal Q}}\}\leq O(\exp(-\Theta(n)))$
which holds by (\ref{eq:Qprob}) and (\ref{eq:AA}).
Thus to show
(\ref{eq:st3}) it suffices to prove
\be{eq:Dconditioned}
\pr\{{\cal D}(t_{m_0})\cap {\cal A}\cap {\cal Q}\}\leq O(\exp(-\Theta(n))).
\ee
We denote the event $({\cal D}(t_{m_0})\cap {\cal A} \cap {\cal Q})$ by ${\cal D}_c(t_{m_0})$.
We first show that given ${\cal D}_c(t_{m_0})$,
there exists with probability one a time instance $t_m$ with  $s_r\leq t_m\leq t_{m_0}$, such that the following events occur:
\be{eq:st4}
{\cal F}(t_{m}) \equiv \left\{ \mb{Q}_{i,1}(t_m)\geq \delta C n \right\}
\ee
and
\be{eq:st5}
{\cal G}(t_{m}) = \left\{ \mb{D}_{i,1}(t_m,t_{m+1})
\leq \mu_{i,1}\bT_{i,1}(t_m,t_{m+1})-{2\delta^2 C n\over\theta} \right\}.
\ee
That is, we claim
\be{eq:FG|D}
\pr\left\{\cup_{\{m:s_r\leq t_m\leq t_{m_0}\}}({\cal F}(t_{m})\cap {\cal G}(t_{m})) \; | \; {\cal D}_c(t_{m_0}) \right\}=1.
\ee
On the other hand we also claim that for each $t_m \in [s_r,s_{r+1}]$
\be{eq:G|F}
\mathbb{P}\left \{{\cal G}(t_{m}) \; | \; {\cal F}(t_{m}) \right\} \le O(\exp(-\Theta(n))),
\ee
which implies that $\mathbb{P} ( {\cal F}(t_{m}) \cap {\cal G}(t_{m})   ) \le O(\exp(-\Theta(n)))$.
Together with (\ref{eq:FG|D}) this would imply
\begin{eqnarray*}
\mathbb{P} \left\{ {\cal D}_c(t_{m_0})   \right\} & \leq &
{\mathbb{P} \left\{ \cup_{s_r\leq t_m\leq t_{m_0}}{\cal F}(t_{m}) \cap {\cal G}(t_{m})\right\}\over
\mathbb{P} \left\{ \cup_{s_r\leq t_m\leq t_{m_0}}{\cal F}(t_{m}) \cap {\cal G}(t_{m})| {\cal D}_c(t_{m_0})  \right\}} \\
& \leq & \sum_{t_m\in[s_r,s_{r+1}]}{\mathbb{P} \left\{{\cal F}(t_{m}) \cap {\cal G}(t_{m})\right\}\over
\mathbb{P} \left\{ \cup_{s_r\leq t_m\leq t_{m_0}}{\cal F}(t_{m}) \cap {\cal G}(t_{m})| {\cal D}_c(t_{m_0})  \right\}} \\
& \leq & O(\exp(-\Theta(n))),
\end{eqnarray*}
where again for the last inequality we use the fact that the number of $t_m$ in the
interval $[s_r,s_{r+1}]$ is at most $\theta/\delta$ and
$\theta/\delta \exp(-\Theta(n))=\exp(-\Theta(n))$.
We could then conclude that (\ref{eq:Dconditioned}) holds and we would be done.
 Thus we need to show (\ref{eq:FG|D}) and (\ref{eq:G|F}).
We start by proving (\ref{eq:G|F}).
Note that during
the time interval $I_m$, policy ${\cal U}(\delta)$ either allocates at least
$\bT_{i,1}(t_m,t_{m+1})$ time units to process class $(i,1)$ jobs, or
all the $\mb{Q}_{i,1}(t_m)>\delta Cn$ jobs initially present are processed. In the second
case  ${\cal G}(t_m)$ does not hold since $\delta Cn>\mu_{i,1}\bT_{i,1}(t_m,t_{m+1})$. In the first case, if
$\bT_{i,1}(t_m,t_{m+1})<{2\delta^2 C n\over\mu_{i,1}\theta}$, then ${\cal G}(t_m)$
obviously does not hold,
since the right-hand side in the inequality in (\ref{eq:st5}) is negative.
Otherwise $\bT_{i,1}(t_m,t_{m+1})\geq \Theta(n)$.  In this case, we can
apply the large deviations
bound (\ref{eq:LDBrate}) which holds by Assumption B. Setting
$\epsilon = 2\delta^2 C n/(\bT_{i,1}(t_m,t_{m+1})\theta)\geq 2\delta C/\theta$
%(\textbf{Is this what you did? Otherwise, I am not sure where the $\theta$ comes from.})
in the bound we obtain:
\[
\mathbb{P} \left \{ \mb{D}_{i,1}(t_m,t_{m+1}) < \mu_{i,1}\bT_{i,1}(t_m,t_{m+1})-
{2\delta^2 C n\over\theta}\right \} \le O(\exp(-\Theta(\bT_{i,1}(t_m,t_{m+1}))))=
O(\exp(-\Theta(n)))
\]
where in the last equation we use $\bT_{i,1}(t_m,t_{m+1})\geq \Theta(n)$ and
as usual, $\delta,C$ and $\theta$ are hidden in the $\Theta(\cdot)$ notation.
%Note that this bound holds for $\bT_{i,1}(t_m,t_{m+1})$ sufficiently large
%(\textbf{unless the large deviations assumption can be restated to hold
%for all $t$)}. We can insure that $\bT_{i,1}(t_m,t_{m+1})$ is sufficiently
%large by increasing the scaling constant $n$. (\textbf{Are we sure this is
%true?})
%Using $\bT_{i,1}(t_m,t_{m+1})\leq t_{m+1}-t_m=\delta n$ and $\mu_{i,1}\leq C$,
%we obtain:
%\begin{equation} \label{eqn:depart}
%\mathbb{P} \left \{ \mb{D}_{i,1}(t_m,t_{m+1}) \ge \mu_{i,1}\bT_{i,1}(t_m,t_{m+1})-
%{\delta^2 Cn\over\theta}
%\right \} \ge 1-\exp(-\Theta(\bT_{i,1}(t_m,t_{m+1}))).
%\end{equation}
%Next, since
%$$ \mu_{i,1}\bT_{i,1}(t_m,t_{m+1})-
%{\delta^2 Cn\over\theta} \ge
%\mu_{i,1}\bT_{i,1}(t_m,t_{m+1})-
%{2\delta^2 Cn\over\theta}, $$
%(\ref{eqn:depart}) and the fact that $\bT_{i,1}(t_m,t_{m+1})=\Theta(n)$
%(recall that the solution $(\bQ(t),\bT(t))$ was obtained by rescaling by $n$)
%yield that $\mathbb{P}(B \; | \; A) \le \exp(-\Theta(n))$ as desired.
We conclude that (\ref{eq:G|F}) holds.

We now prove (\ref{eq:FG|D}). Note that if $t_{m_0}-s_r<\delta n$ then
the right-hand side of the inequality in the event ${\cal D}(t_{m_0})$ is negative and
therefore the events ${\cal D}(t_{m_0})$ and ${\cal D}_c(t_{m_0})$ cannot occur. Thus we assume
there exists at least one $t_m\in [s_r,t_{m_0})$. We have
$\mb{D}_{i,1}(s_r,t_{m_0})\geq \sum_{\{m:s_r\leq t_m\leq t_{m_0-1}\}} \mb{D}_{i,1}(t_m,t_{m+1})$ and
\[
\bD_{i,1}(s_r,t_{m_0})\leq \sum_{\{m:s_r-\delta n\leq t_m\leq t_{m_0-1}\}}\bD_{i,1}(t_m,t_{m+1})
\leq \sum_{\{m:s_r\leq t_m\leq t_{m_0-1}\}}\bD_{i,1}(t_m,t_{m+1})+\delta Cn.
\]
The event ${\cal D}(t_{m_0})$ implies that
there exists a $t_m \in [s_r,t_{m_0-1}]$ such that
\begin{eqnarray}
\mb{D}_{i,1}(t_m,t_{m+1})& < & \mu_{i,1}\bT_{i,1}(t_m,t_{m+1})-{2\delta C^{r+2}n-\delta Cn\over
\left\lceil{t_{m_0}-s_r\over \delta n}\right\rceil} \notag \\
& \leq & \mu_{i,1}\bT_{i,1}(t_m,t_{m+1})-{2\delta C^{r+2}n-\delta Cn\over (\theta/\delta)+1} \label{eq:stmhat} \\
& \leq & \mu_{i,1}\bT_{i,1}(t_m,t_{m+1})-{2\delta^2 C^{r+2}n-\delta^2 C\over \theta+\delta} \label{eq:stCr2}, \\
& \leq & \mu_{i,1}\bT_{i,1}(t_m,t_{m+1})-{2\delta^2 C n\over \theta} \notag,
\end{eqnarray}
where we used $t_{m_0}-s_r\leq \theta n$ in (\ref{eq:stmhat}) and we use in (\ref{eq:stCr2}) the fact that
$\theta\geq 1$ by (\ref{eq:bound_theta}),
$\delta\leq 1$ by (\ref{eq:delta}) and as a result $\theta+\delta\leq 2\theta$ and
$(2C^{r+2}-C)/2\geq (2C^2-C)/2>2C$.

Among  $t_m\in [s_r,t_{m_0})$ select the largest $m$ such that
$\mb{D}_{i,1}(t_m,t_{m+1})\leq \mu_{i,1}\bT_{i,1}(t_m,t_{m+1})-{2\delta^2 Cn\over \theta}$
and denote it by $\hat m$. By the derivation above the set of such $t_m$ is non-empty. Thus,
\be{eq:hatm}
\mb{D}_{i,1}(t_{\hat m},t_{\hat m+1})\leq \mu_{i,1}\bT_{i,1}(t_{\hat m},t_{\hat m+1})-{2\delta^2 Cn\over \theta}.
\ee
Moreover, if $\hat m<m_{0}-1$, then for all $\hat m<m\leq m_{0}-1$, we have

%(\textbf{Doesn't this inequality go the wrong way, given the definition of $\hat{m}$?}):
\[
\mb{D}_{i,1}(t_m,t_{m+1})\geq \mu_{i,1}\bT_{i,1}(t_m,t_{m+1})-{2\delta^2 Cn\over \theta},
\]
or
\begin{equation} \label{eqn:mhat}
\mb{D}_{i,1}(t_{\hat m+1},t_{m_{0}})\geq \mu_{i,1}\bT_{i,1}(t_{\hat m+1},t_{m_{0}})-{({m_{0}}-\hat m)2\delta^2
Cn\over \theta}\geq
\mu_{i,1}\bT_{i,1}(t_{\hat m+1},t_{m_{0}})-2\delta Cn,
\end{equation}
where  $m_{0}-\hat m\leq \theta/\delta$ is used. Note, that the bound (\ref{eqn:mhat})
holds trivially if  $\hat m=m_0-1$.  Next, note that the event ${\cal D}
(t_{m_0})$ jointly with (\ref{eqn:mhat})
implies
\begin{eqnarray} \label{eq:stmhat1}
\mb{D}_{i,1}(s_{r},t_{\hat m+1}) & = & \mb{D}_{i,1}(s_{r},t_{m_{0}})-\mb{D}_{i,1}(t_{\hat m+1},t_{m_{0}})
\nonumber \\
& \leq &
\mu_{i,1}\bT_{i,1}(s_{r},t_{m_{0}})-2\delta C^{r+2}n-(\mu_{i,1}\bT_{i,1}(t_{\hat m+1},t_{m_{0}})-2\delta Cn)
 \\
& = &
\mu_{i,1}\bT_{i,1}(s_{r},t_{\hat m+1})-2\delta C^{r+2}n+2\delta Cn. \nonumber
\end{eqnarray}
Thus
\begin{eqnarray}
\mb{D}_{i,1}(s_{r},t_{\hat m}) & \leq & \mb{D}_{i,1}(s_{r},t_{\hat m+1}) \nonumber \\
& \leq & \mu_{i,1}\bT_{i,1}(s_{r},t_{\hat m+1})-2\delta C^{r+2}n+2\delta Cn  \nonumber\\
& \leq & \mu_{i,1}\bT_{i,1}(s_{r},t_{\hat m})-2\delta C^{r+2}n+3\delta C n, \label{eq:st6}
\end{eqnarray}
where $\mu_{i,1}\bT_{i,1}(t_{\hat m},t_{\hat m+1})\leq \mu_{i,1}\delta n<\delta Cn$ is used.
Now recall from (\ref{eq:stQueueStage1}) that
\begin{equation} \label{eqn:stochdyn}
\mb{Q}_{i,1}(t_{\hat m})=\mb{Q}_{i,1}(s_{r})+\mb{A}_i(s_{r},t_{\hat m})-\mb{D}_{i,1}(s_{r},t_{\hat m}).
\end{equation}
Then conditioned on ${\cal D}_c(t_{m_0})={\cal D}(t_{m_0})\cap {\cal A} \cap {\cal Q}$
and using (\ref{eq:st6}), we obtain
\begin{eqnarray*}
\mb{Q}_{i,1}(t_{\hat m})  \ge \bQ_{i,1}(s_{r})+\lambda_i(t_{\hat m}-s_{r})-\mu_{i,1}\bT_{i,1}(s_{r},t_{\hat m})-
\delta C^{r+2}n-2\delta C n-\delta C n+2\delta C^{r+2}n-3\delta C n.
\end{eqnarray*}
Recall from (\ref{eq:flStage1}) that $\bQ_{i,1}(s_{r})+\lambda_i(t_{\hat m}-s_{r})-\mu_{i,1}\bT_{i,1}(s_{r},t_{\hat m})=\bQ_{i,1}(t_{\hat m})\geq 0$.
Then
\be{eq:st8}
\mb{Q}_{i,1}(t_{\hat m})  \ge
\bQ_{i,1}(t_{\hat m})+\delta C^{r+2}n-6\delta Cn
 >  \delta C n,
\ee
We have established that if the event ${\cal D}_c(t_{m_0})$ holds then (\ref{eq:hatm}) and
(\ref{eq:st8}) hold for some $t_{\hat m}\leq t_{m_0}$. In other words, (\ref{eq:FG|D}) holds.
This completes the proof of (\ref{eq:st3}).

\textbf{Part II.} We now prove  (\ref{eq:st3time}). Fix a $t_m\in [s_r,s_{r+1}]$. Note that the bound (\ref{eq:st3time})
is trivial if
$\bT_{i,1}(s_r,t_m)<3\mu_{\rm max}\delta C^{r+2}n$. So, suppose the previous inequality
does not hold. Let
\be{eq:Delta}
\Delta=\mu_{i,1}\bT_{i,1}(s_{r},t_{m})-2\delta C^{r+2} n\geq \delta C^{r+2}n\geq \Theta(n),
\ee
where we use $\mu_{i,1}\mu_{\rm max}\geq 1$.
Let
\[
\epsilon={\delta Cn\over \Delta}\geq {\delta \over \theta},
\]
where we use $\Delta<\mu_{i,1}\bT_{i,1}(s_{r},t_{m})<C\theta n$.
We condition on the event ${\cal D}(t_m)$, which by (\ref{eq:st3}) holds with probability at least $1-\exp(-\Theta(n))$
and use large deviations Assumption B with the $\epsilon$ above to obtain
\begin{eqnarray*}
\mathbb{P} \left\{ \mb{T}_{i,1}(s_{r},t_m)\geq \mu_{i,1}^{-1}\Big(\mu_{i,1}\bT_{i,1}(s_{r},t_{m})-
2\delta C^{r+2} n\Big)-\delta C n \; \Big| \; {\cal D}(t_m) \right\}  & = & \\
\mathbb{P} \left\{ \mb{T}_{i,1}(s_{r},t_m)\geq \mu_{i,1}^{-1}\Delta-\epsilon\Delta
\; \Big| \; \mb{D}(s_r,t_m)\geq\Delta \right\} &\geq & \\
1-O(\exp(-\Theta(\Delta))) & \geq & \\
1-O(\exp(-\Theta(n))),
\end{eqnarray*}
where the last inequality follows from the last inequality in (\ref{eq:Delta}).
To finish the argument, we observe that $\mu_{i,1}^{-1}2\delta C^{r+2}+\delta C n<3\mu_{\max}\delta C^{r+2}n$. \qed

We now establish a similar lower bound for classes corresponding to stages two and higher.

\begin{lemma}\label{lemma:DepartureStage2+}
For every $i\leq I$, $j\leq J_i$ and $m$ such that $s_r\leq t_m\leq s_{r+1}$
\be{eq:st3stage2+}
\mathbb{P}\{\mb{D}_{i,j}(s_{r},t_{m})\geq \bD_{i,j}(s_{r},t_{m})-2\delta j  C^{r+2} n\}\geq 1-\exp(-\Theta(n)),
\ee
and
\be{eq:st3timeStage2+}
\mathbb{P}\{\mb{T}_{i,j}(s_{r},t_{m})\geq \bT_{i,j}(s_{r},t_{m})-3\mu_{\max}\delta j C^{r+2} n\}
\geq 1-\exp(-\Theta(n)).
\ee
\end{lemma}

\proof The proof is very similar to the one for Lemma \ref{lemma:DepartureStage1}. We only highlight the
differences. The proof is done by induction in $j$; the base case $j=1$ is covered by
Lemma \ref{lemma:DepartureStage1}.
So let us fix a $j > 1$ and
assume that the assertion holds for all $(i,j')$
with $j'\leq j-1$. We again define an event related to the inequality
inside $\mathbb{P}$ in (\ref{eq:st3stage2+}). For a class $(i,j)$ and any
time $t_{m_0}$ with $s_r\leq t_{m_0} \leq s_{r+1}$ let
\be{eq:Dj}
{\cal D}(t_{m_0}) = \{ \mb{D}_{i,j}(s_{r},t_{m})< \bD_{i,j}(s_{r},t_{m})-2\delta j  C^{r+2} n \}.
\ee
We need to show $\pr\left\{{\cal D}(t_{m_0})\right\}\leq O(\exp(-\Theta(n)))$ for every $t_{m_0}\in [s_r,s_{r+1}]$.
As in Lemma \ref{lemma:DepartureStage1} we introduce the event
\be{eq:Qj}
{\cal Q} \equiv \left\{\mb{Q}_{i,j}(s_{r})\geq \bQ_{i,j}(s_{r})-\delta C^{r+2} n-2\delta Cn \right\},
\ee
but instead of the event ${\cal A}$ defined by (\ref{eq:A}), consider
\be{Dj-1}
{\cal D} \equiv \left\{\forall \; t_m\in[s_r,s_{r+1}]: \mb{D}_{i,j-1}(s_r,t_m) \geq \bD_{i,j-1}(s_r,t_m)-2\delta (j-1)C^{r+2} n \right\}.
\ee
Again using Lemma \ref{lemma:sr} (and the ``outer'' inductive assumption)
we obtain  $\pr\left\{{\cal Q} \right\}\geq 1-\exp(-\Theta(n))$ and by the inductive
assumption on $j$, $\pr\left\{{\cal D} \right\}\geq 1-O(\exp(-\Theta(n)))$
(where, as before, we sum several expressions of the order $O(\exp(-\Theta(n)))$
over $t_m\in [s_r,s_{r+1}]$ to get again $O(\exp(-\Theta(n)))$).
Next, let ${\cal D}_c(t_{m_0})={\cal D}(t_{m_0})\cap {\cal Q} \cap {\cal D}$.
We need to show $\pr\left\{{\cal D}_c(t_{m_0})\right\}\leq O(\exp(-\Theta(n)))$.
For every $t_m\in [s_r,t_{m_0}]$, we introduce the event  ${\cal F}(t_m)$
as in (\ref{eq:st4}), except
$\mb{Q}_{i,j}$ is used instead of $\mb{Q}_{i,1}$. Finally, we introduce ${\cal G}(t_m)$, defined
as follows
\[
\ {\cal G}(t_{m}) = \left\{ \mb{D}_{i,j}(t_m,t_{m+1})
\leq \mu_{i,j}\bT_{i,j}(t_m,t_{m+1})-{2\delta^2 C nj\over\theta} \right\}.
\]
Arguing as in the proof of Lemma \ref{lemma:DepartureStage1},
we claim that (\ref{eq:FG|D}) and (\ref{eq:G|F}) hold with the new event definitions.
The proof of (\ref{eq:G|F}) is identical to the one of Lemma \ref{lemma:DepartureStage1}. For (\ref{eq:FG|D})
we repeat the argument until we get to (\ref{eq:st6}), instead of which we get
\begin{eqnarray} \label{eq:st66}
\mb{D}_{i,j}(s_{r},t_{\hat m}) \leq
\mu_{i,1}\bT_{i,j}(s_{r},t_{\hat m})-2j\delta C^{r+2}n+3\delta C n.
\end{eqnarray}

Then we obtain
\begin{eqnarray}
\mb{Q}_{i,j}(t_{\hat m}) &= & \mb{Q}_{i,j}(s_{r})+\mb{D}_{i,j-1}(s_{r},t_{\hat m})-\mb{D}_{i,j}(s_{r},t_{\hat m}) \label{eq:Qijst1} \\
& \geq & \bQ_{i,j}(s_{r})+\mu_{i,j-1}\bT_{i,j-1}(s_{r},t_{\hat m})-\mu_{i,j}\bT_{i,j}(s_{r},t_{\hat m}) \notag \\
& - & \delta C^{r+2} n -2\delta Cn-2\delta (j-1)C^{r+2} n+2\delta j C^{r+2} n-3\delta C n \label{eq:Qijst2} \\
& = & \bQ_{i,j}(t_{\hat m})+\delta C^{r+2} n-5\delta C n \label{eq:Qijst3} \\
& \geq & \delta C n \notag,
\end{eqnarray}
where (\ref{eq:stQueueStagej}) is used for (\ref{eq:Qijst1}), conditioning on
${\cal D}_c(t_{m_0})$ is used in
(\ref{eq:Qijst2}), and (\ref{eq:flStagej}) is used in (\ref{eq:Qijst3}). This proves (\ref{eq:FG|D})
and completes the proof of (\ref{eq:st3stage2+}). The proof of the lower bound for
$\mb{T}_{i,j}(\cdot)$ follows the proof of Lemma \ref{lemma:DepartureStage1},
almost line for line. \qed

\subsubsection{Upper bounds on the departure processes}\label{subsection:DepartureUpper}
In this subsection we obtain upper bounds,
similar to the bounds in Lemmas \ref{lemma:DepartureStage1} and \ref{lemma:DepartureStage2+},
on the cumulative departures
$\mb{D}_{i,j}(s_{r},t_m)$, for values of $m$ such that $s_r<t_m\leq s_{r+1}$.

For every station  $\sigma$, by construction of the sequence
$s_0,s_1,\ldots,s_M$, we have either
$\bQ_{\sigma}(t)>0$ for all $s_{r}<t<s_{r+1}$ or $\bQ_{\sigma}(t)=0$ for all $s_{r}<t<s_{r+1}$.
We consider these cases separately.

\begin{lemma}\label{lemma:DepartureUpper}
Given any station $\sigma$, suppose the interval $(s_r,s_{r+1})$ is such that
$\bQ_{\sigma}(t)>0$ for all $s_{r}<t<s_{r+1}$.
Then, for every class $(i,j)\in\sigma$  and every $t_m\in [s_r,s_{r+1}]$
\be{eq:stUpper22T}
\mathbb{P}\{\mb{T}_{i,j}(s_{r},t_m)\leq \bT_{i,j}(s_{r},t_m)+3\delta |\sigma|\mu_{\max}J_{\max}C^{r+2}n\}\geq 1-O(\exp(-\Theta(n))),
\ee
and
\be{eq:stUpper22D}
\mathbb{P}\{\mb{D}_{i,j}(s_{r},t_m)\leq \bD_{i,j}(s_{r},t_m)+4\delta |\sigma|\mu_{\max}^2J_{\max}C^{r+2}n\} \geq 1-O(\exp(-\Theta(n))).
\ee
\end{lemma}

\proof
Given any station $\sigma$, suppose that $\bQ_{\sigma}(t)>0$ for all $s_{r}<t<s_{r+1}$. By
the non-idling constraint (\ref{eq:flNon-idling}) we have that
\be{eq:stUpper1}
\sum_{(i,j)\in\sigma}\bT_{i,j}(s_{r},s_{r+1})=s_{r+1}-s_{r}.
\ee
Fix any $t_m\in [s_r,s_{r+1}]$ and fix any  class $(i,j)\in\sigma$.
Applying Lemmas \ref{lemma:DepartureStage1} and \ref{lemma:DepartureStage2+} to  $s_r$ and $t_m$
%(\textbf{How do we know this is possible?}),
we have that with probability at least $1-O(\exp(-\Theta(n)))$, for every class $(i',j')\in\sigma$,
$$\mb{T}_{i',j'}(s_{r},t_m)\geq\bT_{i',j'}(s_{r},t_m)-3\delta j' \mu_{\max}C^{r+2} n\geq
\bT_{i',j'}(s_{r},t_m)-3\delta \mu_{\max}J_{\max}C^{r+2} n,$$
where we use (\ref{eq:Jmax}) in the second inequality.
Applying (\ref{eq:stUpper1}) and the feasibility inequality (\ref{eq:stFeasibility}), we obtain
that with probability at least $1-O(\exp(-\Theta(n)))$,
\be{eq:stUpper2}
\mb{T}_{i,j}(s_{r},t_m)\leq t_m-s_{r}-
\sum_{(i,j)\neq (i',j')\in\sigma}\mb{T}_{i',j'}(s_{r},t_m)\leq
\ee
\[
t_m-s_{r}-\sum_{(i,j)\neq (i',j')\in\sigma}(\bT_{i',j'}(s_{r},t_m)-3\delta \mu_{\max} J_{\max}C^{r+2}n)\leq
\bT_{i,j}(s_{r},t_m)+3\delta |\sigma|\mu_{\max}J_{\max}C^{r+2}n.
\]
%\textbf{Did we define the notation $|\sigma|$ somewhere?}
Let us define the event ${\cal T}$ as follows:
\[
{\cal T}=\left\{\mb{T}_{i,j}(s_{r},t_m)\leq \bT_{i,j}(s_{r},t_m)+
3\delta |\sigma|\mu_{\max}J_{\max}C^{r+2}n\right\},
\]
that is $\pr\left\{{\cal T} \right\}\geq 1-O(\exp(-\Theta(n)))$ per (\ref{eq:stUpper2}).
Next, consider
\begin{eqnarray*}
\pr\left\{\mb{D}_{i,j}(s_{r},t_m)>\mu_{i,j}\bT_{i,j}(s_{r},t_m)+4\delta
|\sigma|\mu^2_{\max}J_{\max}C^{r+2}n \; | \; {\cal T} \right\} &\leq & \\
\pr\left\{\mb{D}_{i,j}(s_{r},t_m)>\mu_{i,j}(\bT_{i,j}(s_{r},t_m)+3\delta
|\sigma|\mu_{\max}J_{\max}C^{r+2}n)+
\delta |\sigma|\mu^2_{\max}J_{\max}C^{r+2}n \; | \; {\cal T} \right\}. & &
\end{eqnarray*}
Applying the large deviation bound (\ref{eq:LDBrate}) with
\[
\epsilon={\delta |\sigma|\mu^2_{\max}J_{\max}C^{r+2}n\over \bT_{i,j}(s_{r},t_m)+3\delta |\sigma|\mu_{\max}J_{\max}C^{r+2}n}
\geq {\delta |\sigma|\mu^2_{\max}J_{\max}C^{r+2} \over \theta+3\delta |\sigma|\mu_{\max}J_{\max}C^{r+2}},
\]
(where $\bT_{i,j}(s_{r},t_m)\leq \theta n$ is used), we obtain
\begin{eqnarray*}
\pr\left\{\mb{D}_{i,j}(s_{r},t_m)>\mu_{i,j}\bT_{i,j}(s_{r},t_m)+4\delta |\sigma|\mu^2_{\max}J_{\max}C^{r+2}n\right|T\} &\leq & \\
\exp\Big(-\Theta\left(\bT_{i,j}(s_{r},t_m)+3\delta |
\sigma|\mu_{\max}J_{\max}C^{r+2}n\right)\Big)\leq \exp(-\Theta(n)),
\end{eqnarray*}
where the constant $\epsilon$ is hidden in $\Theta(\cdot)$.
Taking this together with $\pr\left\{{\cal T} \right\}\geq 1-O(\exp(-\Theta(n)))$,
we have proven the lemma.
\qed

We now analyze stations $\sigma$ for which the fluid amount stays zero
during the interval $[s_r,s_{r+1}]$. In the
following lemma we obtain an analogue of Lemma \ref{lemma:DepartureUpper} for
this second case.

\begin{lemma}\label{lemma:DepartureUpperCaseII}
Given any station $\sigma$, suppose the interval $(s_r,s_{r+1})$ is such that $\bQ_{\sigma}(t)=0,$
for all $s_{r}<t<s_{r+1}$. Then for every $(i,j)\in\sigma$ and every $t_m\in [s_r,s_{r+1}]$
\be{eq:stUpper2CaseII}
\mathbb{P}\{\mb{D}_{i,j}(s_{r},t_m)\leq \bD_{i,j}(s_{r},t_m)+5\delta  j\mu_{\max}^2IJ_{\max}^2C^{r+2}n\}\geq 1-O(\exp(-\Theta(n))).
\ee
\end{lemma}

\proof

Consider any station $\sigma$ such that $\bQ_{\sigma}(t)=0$, for all $s_{r}<t<s_{r+1}$.
Applying fluid equations (\ref{eq:flStage1}) and (\ref{eq:flStagej})
we obtain that for every class $(i,j)\in\sigma$ and every $t_m\in [s_r,s_{r+1}]$
\be{eq:stUpper1CaseII}
\bD_{i,j-1}(s_{r},t_m)=\bD_{i,j}(s_{r},t_m),
\ee
where for the case $j=1$, $\bD_{i,j-1}(\cdot)$ is understood as $\bA_i(\cdot)$.

The proof now proceeds by induction in $j$.
We start with the base step, $j=1$.
So, consider any class $(i,1)\in\sigma$. Applying (\ref{eq:stUpper1CaseII}), we have
\be{eq:stUpper3CaseII}
\bD_{i,1}(s_{r},t_m)=\bA_i(s_{r},t_m).
\ee
Applying (\ref{eq:stQueueStage1}) we have
\be{eq:stUpper4CaseII}
\mb{Q}_{i,1}(t_m)=\mb{Q}_{i,1}(s_{r})+\mb{A}_i(s_{r},t_m)-\mb{D}_{i,1}(s_{r},t_m)\geq 0.
\ee
Next, let $t_{m'} \equiv \max \{ m: t_m \leq s_r \}$.
In particular, $0\leq s_r-t_{m'}\leq \delta n$. We then have
\begin{eqnarray}
\mathbb{P} \{ \mb{A}_i(s_{r},t_m)>\bA_i(s_{r},t_m)+2\delta Cn \} & \leq & \mathbb{P} \{ \mb{A}_i(t_{m'},t_m)>\bA_i(s_{r},t_m)+2\delta Cn \} \\
&\leq & \mathbb{P} \{ \mb{A}_i(t_{m'},t_m)>\bA_i(t_{m'}+\delta n,t_m)+2\delta Cn \} \\
&\leq & \mathbb{P} \{ \mb{A}_i(t_{m'},t_m)>\bA_i(t_{m'},t_m)+\delta Cn \} \label{eq:srtm'}, \label{eq:Asrtm}
\end{eqnarray}
where $\bA_i(t_{m'},t_{m'}+\delta n)\leq \lambda_i\delta n<C\delta n$ is used in (\ref{eq:srtm'}).
Note that when $t_{m'}=t_m$ the probability in (\ref{eq:srtm'}) is zero since the left-hand side is negative.
Thus we assume $t_{m'}>t_m$. Let $\epsilon=\delta C n/(t_m-t_{m'})\geq \delta Cn/(\theta n)=\delta C/\theta>0$. Using the large
deviations Assumption B with this $\epsilon$ we obtain
\begin{eqnarray*}
\mathbb{P} \left\{ \mb{A}_i(t_{m'},t_m)>\bA_i(t_{m'},t_m)+\delta Cn \right\}  & = &
\mathbb{P} \{ \mb{A}_i(t_{m'},t_m)>\bA_i(t_{m'},t_m)+\epsilon (t_m-t_{m'}) \} \\
& \leq & O(\exp(-\Theta(t_m-t_{m'}))) \\
& \leq & O(\exp(-\Theta(n))),
\end{eqnarray*}
where $t_m-t_{m'}\geq \delta n$ is used in the last inequality. Combining this bound with (\ref{eq:Asrtm}) we obtain
\be{eq:Asrtm1}
\mathbb{P} \{ \mb{A}_i(s_{r},t_m)>\bA_i(s_{r},t_m)+2\delta Cn \}  \leq O(\exp(-\Theta(n))).
\ee
%\textbf{Again, for $s_{r+1}-s_r$ sufficiently large. Also, are we setting $\epsilon = C$ here?}
Applying Lemma \ref{lemma:sr} we have
$$ \mathbb{P} \{\mb{Q}_{i,1}(s_{r})\leq \bQ_{i,1}(s_{r})+\delta C^{r+2}n+2\delta Cn \} \geq 1-O(\exp(-\Theta(n))).$$
By our assumption that $\bQ_{\sigma}(t)=0$ for all $s_r<t<s_{r+1}$ and by continuity, we have
$\bQ_{\sigma}(s_{r})=0$. Using this fact, we now have
\be{eq:Qi1}
\mathbb{P} \{\mb{Q}_{i,1}(s_{r})\leq \delta C^{r+2}n+2\delta C n \} \geq 1-O(\exp(-\Theta(n))).
\ee
Now, from (\ref{eq:stUpper4CaseII}), we have
\begin{eqnarray*}
\mb{D}_{i,1}(s_{r},t_m) &  = & \mb{Q}_{i,1}(s_{r})+\mb{A}_i(s_{r},t_m) - \mb{Q}_{i,1}(t_m) \\
& \leq & \mb{Q}_{i,1}(s_{r})+\mb{A}_i(s_{r},t_m).
\end{eqnarray*}
Applying (\ref{eq:Asrtm1}) and (\ref{eq:Qi1}) we obtain
that with probability at least $1-O(\exp(-\Theta(n)))$,
\[
\mb{D}_{i,1}(s_{r},t_m)\leq \bA_i(s_{r},t_m)+2\delta Cn+\delta C^{r+2}n+2\delta Cn\leq \bA_i(s_{r},t_m)+5\delta\mu_{\max}^2IJ_{\max}^2C^{r+2}n.
\]
Combining this with (\ref{eq:stUpper3CaseII}), we obtain the required bound.
This completes the proof of the
base step.

We now prove the inductive step. So, fix $j > 1$ and suppose that the assertion
holds for $1, 2, \ldots, j-1$. We now consider a particular class $(i,j)\in\sigma$.
We have from (\ref{eq:stQueueStagej})
\be{eq:stUpper4CaseII0}
\mb{Q}_{i,j}(t_m)=\mb{Q}_{i,j}(s_{r})+\mb{D}_{i,j-1}(s_{r},t_m)-\mb{D}_{i,j}(s_{r},t_m)\geq 0.
\ee

Again, by Lemma \ref{lemma:sr} we have
$$ \mathbb{P} \{ \mb{Q}_{i,j}(s_{r})\leq \bQ_{i,j}(s_{r})+\delta C^{r+2}n+2\delta C n \}
\geq 1-O(\exp(-\Theta(n))),$$
which implies
\be{eq:Qij2}
 \mathbb{P} \{\mb{Q}_{i,j}(s_{r})\leq \delta C^{r+2}n+2\delta C n \} \geq 1-O(\exp(-\Theta(n))),
\ee
again since $\bQ_{\sigma}(s_{r})=0$.

Consider the station $\sigma_{\nu'}$ containing $(i,j-1)$. If $\sigma_{\nu'}$ is also such that $\bQ_{\sigma'}(t)=0$ for
all $s_{r}<t<s_{r+1}$ (for example when $\nu'=\nu$), then by the inductive assumption on $j$,
$$ \mathbb{P} \{ \mb{D}_{i,j-1}(s_{r},t_m)\leq
\bD_{i,j-1}(s_{r},t_m)+5\delta (j-1)\mu_{\max}^2IJ_{\max}^2C^{r+2}n \}
\geq 1-O(\exp(-\Theta(n))). $$
Otherwise, $\sigma_{\nu'}$ is such that $\bQ_{\sigma_{\nu'}}(t)>0$ for all $s_{r}<t<s_{r+1}$.
Then Lemma \ref{lemma:DepartureUpper} becomes applicable, and applying  (\ref{eq:stUpper22D}) to $\sigma_{\nu'}$, we
have that with probability at least $1-O(\exp(-\Theta(n)))$,
\begin{eqnarray*}
 \mb{D}_{i,j-1}(s_{r},t_m)& \leq & \bD_{i,j-1}(s_{r},t_m)+
4\delta |\sigma_{\nu'}|\mu_{\max}^2J_{\max}C^{r+2}n \\
& < & \bD_{i,j-1}(s_{r},t_m)+ 5(j-1)\delta \mu_{\max}^2IJ_{\max}^2C^{r+2}n,
\end{eqnarray*}
where we use $|\sigma_{\nu'}|\leq \sum_iJ_i\leq IJ_{\max}$ and $j > 1$.
Hence, in either case we have a probabilistic bound on $\mb{D}_{i,j-1}(s_{r},t_m)$.
Combining this with
(\ref{eq:stUpper4CaseII0}) and (\ref{eq:Qij2}) we obtain that with probability at least $1-O(\exp(-\Theta(n)))$,
\begin{eqnarray*}
\mb{D}_{i,j}(s_{r},t_m) & \leq & \bD_{i,j-1}(s_{r},t_m)+\delta C^{r+2}n+2\delta C n+5\delta (j-1)\mu_{\max}^2IJ^2_{\max}C^{r+2}n \\
& < &
\bD_{i,j-1}(s_{r},t_m)+\delta(3+5 (j-1)\mu_{\max}^2IJ^2_{\max})C^{r+2}n \\
& < &
\bD_{i,j-1}(s_{r},t_m)+5\delta j\mu_{\max}^2IJ_{\max}^2C^{r+2}n.
\end{eqnarray*}
Finally, recalling (\ref{eq:stUpper1CaseII}) we obtain the desired
bound. This completes the proof of the inductive step. \qed

\vspace{.1in}

With the lemmas above in hand,
we are now ready to finish the proof of Proposition \ref{prop:StochasticClose},
by completing the outer inductive step on $r$.

\vspace{0.1in}

\proof[Proof of Proposition \ref{prop:StochasticClose}]

Fix any $t_m\in [s_r,s_{r+1}]$.
By Lemma \ref{lemma:sr}
$$ \mathbb{P} \{
|\mb{Q}_{i,j}(s_{r})-\bQ_{i,j}(s_{r})| \leq \delta C^{r+2}n+2\delta C n \}
\geq 1-O(\exp(-\Theta(n))). $$
Next, for any class $(i,j)$ recall that we have
\begin{equation} \label{eqn:finalst}
\mb{Q}_{i,j}(t_m)=\mb{Q}_{i,j}(s_r)+\mb{D}_{i,j-1}(s_r,t_m)-\mb{D}_{i,j}(s_r,t_m),
\end{equation}
with $\mb{D}_{i,j-1}(\cdot)$ replaced by $\mb{A}_i(\cdot)$ when $j=1$.  Combining Lemmas
\ref{lemma:DepartureStage1},\ref{lemma:DepartureStage2+},\ref{lemma:DepartureUpper}, and \ref{lemma:DepartureUpperCaseII},
we obtain that,
$$ \mathbb{P} \{ |\mb{D}_{i,j}(s_{r},t_m)-\bD_{i,j}(s_{r},t_m)|\leq 5j\delta \mu_{\max}^2IJ_{\max}^2C^{r+2}n \}
\geq 1-O(\exp(-\Theta(n))). $$
Recalling (\ref{eq:Alower}) and (\ref{eq:Asrtm1}) we have that
$$ \mathbb{P} \{ |\mb{A}_i(s_{r},t_m)-\bA_i(s_{r},t_m)|\leq 2\delta Cn \}
\geq 1-O(\exp(-\Theta(n))). $$

Combining the previous two bounds with (\ref{eqn:finalst}) and the fluid analogs
(\ref{eq:flStage1}) and (\ref{eq:flStagej})
we obtain that, with probability at least $1-O(\exp(-\Theta(n)))$,
\[
|\mb{Q}_{i,j}(t_m)-\bQ_{i,j}(t_m)|\leq \delta C^{r+2}n+2\delta C n+10\delta j \mu_{\max}^2IJ_{\max}^2C^{r+2}n
<\delta C^{r+3}n,
\]
where the bound (\ref{eq:C}) is used. This completes the proof of
Proposition \ref{prop:StochasticClose}. \qed

\subsection{Proof of the Main Theorem} \label{sec:mp}

In this subsection we present the final two proofs. The next theorem is
last result needed before proving Theorem~\ref{theorem:Main}.

\begin{theorem}\label{theorem:mainStoch}
Suppose the fluid model of a stochastic multitype network is not
globally weakly stable and satisfies the FDP property. Then for any initial
state $\mb{Q}(0)=(q,z_1,z_2)\in
\mathbb{Z}_+^d\times\Re_+^{I+d}$, under the non-idling scheduling policy ${\cal U}$ we have
\be{eq:st1}
\mathbb{P} \left\{\big|\big|\mb{Q}(\theta ||q||)\big|\big|\geq 2||q|| \right\} \ge
1-O(\exp(-\Theta(||q||))), \ee and \be{eq:st11} \mathbb{P}
\left\{\inf_{0\leq t\leq\theta ||q||}||\mb{Q}(t)||\geq {||q||\over 4}
\max\left({\gamma\over C},1\right) \right\} \ge
1-O(\exp(-\Theta(||q||))). \ee
\end{theorem}

\proof
We first prove (\ref{eq:st11}). Fix any $t_m$ and find the $s_r$
such that  $t_m\in [s_r,s_{r+1}]$.
We have
\begin{eqnarray*}
\mathbb{P}\left\{\sup_{t_m\leq t\leq t_{m+1}}||\mb{Q}(t)-\bQ(t)|| > 3C^{r+4}\delta n\right\} & \leq &
\mathbb{P}\left\{\sup_{t_m\leq t\leq t_{m+1}}||\mb{Q}(t)-\mb{Q}(t_m)|| > C\delta n\right\} \\
& + &\mathbb{P}\left\{||\mb{Q}(t_m)-\bQ(t_m)|| > C^{r+4}\delta n\right\} \\
& + &\mathbb{P}\left\{\sup_{t_m\leq t\leq t_{m+1}}||\bQ(t_m)-\bQ(t)|| > C\delta n\right\} \\
& \leq & O(\exp(-\Theta(n))).
\end{eqnarray*}

In fact, observe that the last probability in the right-hand side above is equal to
zero by (\ref{eq:Cdeltanbar}) of Lemma \ref{lemma:QQt_m}. The first
probability in the right-hand side is at most $O(\exp(-\Theta(n)))$  by (\ref{eq:mbQ2}) of Lemma \ref{lemma:QQt_m},
and the second probability is also at most  $O(\exp(-\Theta(n)))$ by
(\ref{eq:st2}) of Proposition~\ref{prop:StochasticClose} and the fact that $\sum J_i<C$.

Combining the inequality above with (\ref{eq:bQn1>3}) we obtain
\[
\mathbb{P}\left\{\inf_{t_m\leq t\leq t_{m+1}}||\mb{Q}(t)|| < {n\over 2}
\min\left({\gamma\over C},1\right)-3C^{r+4}\delta n\right\}\leq O(\exp(-\Theta(n))).
\]
>From (\ref{eq:delta}) and since $r\leq M-1$ we have
\[
{n\over 2}\min\left({\gamma\over C},1\right)-3C^{r+4}\delta n\geq {n\over 4}\min
\left({\gamma\over C},1\right).
\]

Thus
\[
\mathbb{P}\left\{\inf_{t_m\leq t\leq t_{m+1}}||\mb{Q}(t)|| < {n\over 4}\min({\gamma\over C},1)\right\}\leq O(\exp(-\Theta(n))).
\]
By summing over all $m=0,1,\ldots,\lceil\theta/\delta\rceil$ we obtain
\[
\mathbb{P}\left\{\inf_{0\leq t\leq \theta n}||\mb{Q}(t)|| <
{n\over 4}\min\left({\gamma\over C},1\right)\right\}\leq \left\lceil{\theta\over\delta}\right\rceil O(\exp(-\Theta(n)))=
O(\exp(-\Theta(n))),
\]
where the last equality follows since by (\ref{eq:bound_theta}) and (\ref{eq:delta}), the value of $\lceil\theta/\delta\rceil$ is
bounded above by a constant. Recall, finally, that $||q||=n$.
This completes the proof of (\ref{eq:st11}).

We now prove (\ref{eq:st1}). Find the largest $t_m\leq \theta n$. In particular $\theta n -t_m\leq \delta n$.
Applying (\ref{eq:mbQ2}) with $t=\theta n$, we obtain
\be{eq:mbQtheta}
\mathbb{P}\left\{||\mb{Q}(\theta n)-\mb{Q}(t_m)| > C\delta n\right\} \leq O(\exp(-\Theta(n))).
\ee
Applying (\ref{eq:Cdeltanbar}) at $t=\theta n$ we obtain
\be{eq:Cdeltanbar1}
||\bQ(\theta n)-\bQ(t_m)||\leq C\delta n.
\ee
Applying (\ref{eq:st2}) to the $t_m$ chosen above we obtain
\be{eqn:qij}
\mathbb{P}\{|\mb{Q}_{i,j}(t_m)-\bQ_{i,j}(t_m)|> \delta C^{M+2} n\}\leq O(\exp(-\Theta(n))).
\ee
Next, we note that
\begin{eqnarray*}
\mathbb{P}\{||\mb{Q}(t_m)-\bQ(t_m)||>\delta C^{M+3} n\} & \leq &
\mathbb{P}\left\{ \bigcup_{i,j} \quad |\mb{Q}_{i,j}(t_m)-\bQ_{i,j}(t_m)| > \frac{\delta C^{M+3} n}{IJ_{\max}}\right\} \\
& \leq & \mathbb{P}\left\{ \bigcup_{i,j} \quad |\mb{Q}_{i,j}(t_m)-\bQ_{i,j}(t_m)| > \frac{\delta C^{M+3} n}{C}\right\} \\
& = & \mathbb{P}\left\{ \bigcup_{i,j} \quad |\mb{Q}_{i,j}(t_m)-\bQ_{i,j}(t_m)| > \delta C^{M+2} n\right\} \\
& \leq & \sum_{i,j} \mathbb{P}\left\{|\mb{Q}_{i,j}(t_m)-\bQ_{i,j}(t_m)| > \delta C^{M+2} n\right\} \\
& \leq & O(\exp(-\Theta(n))).
\end{eqnarray*}
In the last step, we employ (\ref{eqn:qij}) and then sum over all $i$ and $j$
to obtain a new exponential bound. Combining (\ref{eq:mbQtheta}), (\ref{eq:Cdeltanbar1}) and the last
bound, we obtain
\[
\mathbb{P}\{||\mb{Q}(\theta n)-\bQ(\theta n)||>3\delta C^{M+3} n\}\leq O(\exp(-\Theta(n))).
\]
Since $\delta <1/(3C^{M+3})$
we obtain:
\[
\mathbb{P}\{||\mb{Q}(n\theta)-\bQ(n\theta)||> n\} \leq O(\exp(-\Theta(n))).
\]
Recalling from (\ref{eq:bQn>3}) that $||\bQ(n\theta)||\geq 3n$ and recalling $||q||=n$, we obtain
\[
\mathbb{P}\{||\mb{Q}(n\theta)||< 2||q||\} \leq O(\exp(-\Theta(n))),
\]
which implies (\ref{eq:st1}). This completes the proof of Theorem \ref{theorem:mainStoch}. \qed

It should be noted that the
constant ``2'' which appears in (\ref{eq:st1}) is completely arbitrary. In all of the
proofs in which the constant appears, it can be replaced
by any constant greater than unity. We are now ready to prove our main result.

\proof[Proof of Theorem \ref{theorem:Main}]
We fix a large value $n_0$ (the actual value will be specified later).
Consider any initial state $(q,z_1,z_2)$ with $||q||\geq n_0$.
We  apply the policy ${\cal U}$
for the time interval $[0,\theta_0]$ where $\theta_0=\theta ||q||$. If at time $\theta_0 $ the resulting
state $\mQ(\theta_0)$ is such that $||\mQ(\theta_0)||\geq 2n_0$, then we apply the
policy ${\cal U}$ again with $q$ reset to $\mQ(\theta_0)$, till the corresponding time $\theta_1=\theta_0+\theta ||\mQ(\theta_0)||$.
If again $||\mQ(\theta_1)||\geq 2||\mQ(\theta_0)||\geq 4n_0$,
we continue with policy ${\cal U}$ until the corresponding time $\theta_2$ and check
whether $||\mQ(\theta_2)||\geq 2||\mQ(\theta_1)||\geq 8n_0$, and so on. Either this process
continues indefinitely or for some time instance $\theta_i$ we get
$||\mQ(\theta_i)||< 2||\mQ(\theta_{i-1})||$. Set $\theta_{-1}=0$ by convention.
Let ${\cal E}_m, m=0,1,\ldots,$ denote the event
$||\mQ(\theta_{i})||\geq 2||\mQ(\theta_{i-1})||$ and
\be{eq:Qinfty}
\inf_{\theta_{i-1}\leq t\leq \theta_{i}}||\mQ(t)||\geq {||\mQ(\theta_{i-1})||\over 4}\max\left({\gamma\over C},1\right),
\ee
for all $i\leq m$. In particular, the event implies $||\mQ(\theta_{m})||\geq 2^{m+1}n_0\geq n_0$.
Let ${\cal E}^1=\cap_m {\cal E}_m$,
that is ${\cal E}^1$ implies  that the process of exceeding the bounds continues indefinitely.
We now show that $\pr\left\{{\cal E}^1\right\}\geq\alpha>0$, where
$\alpha$ depends only on the parameters of the model and on $n_0$ and $\gamma$
(and is independent for example from the components $z_1,z_2$).
%It remains to show that $\pr\{{\cal E}\}\geq\alpha>0$ for some $\alpha$ which depends only on $n_0,\gamma$ and the parameters of the model,
%and to show that ${\cal E}$ implies (\ref{eq:rate}).
By (\ref{eq:st1}) and (\ref{eq:st11}) of Theorem \ref{theorem:mainStoch},
the probability of the event ${\cal E}^1 \equiv \cap_{m=0}^{\infty}{\cal E}_i$
is at least
\begin{eqnarray*}
1-2\sum_{m=0}^{\infty}O(\exp(-\Theta(2^m||q||))) & > & 1-\sum_{m=0}^{\infty}O(\exp(-\Theta((m+1)||q||))) \\
& = &1-{O(\exp(-\Theta(||q||)))\over 1-O(\exp(-\Theta(||q||)))} \\
& > & 1-O(e^{-\Theta(n_0)}).
\end{eqnarray*}
We take $n_0$ sufficiently large so that $\alpha\equiv 1-O(e^{-\Theta(n_0)})>0$. The parameters hidden in $\Theta(\cdot)$
depend only on the parameters of the model (including the large deviations parameters $V,L$) and $\gamma$.
Thus the probability of ${\cal E}^1$ is positive (and in fact is close to unity), provided
that $n_0$ is sufficiently large.

Next, we show that the event ${\cal E}^1$ implies
\be{eq:rate}
\liminf_t{||\mQ(t)||\over t}\geq {\max\left({\gamma\over C},1\right)\over 8\max(1,{3\over \gamma})}>0.
\ee
We first show by induction in $m$ that
event ${\cal E}^1$ implies   $||\mQ(\theta_m)||/\theta_m\geq 1/\theta$ for all $m\geq 1$.
When $m=1$ the ratio is at least
$2||q||/(\theta ||q||)>1/\theta$.
Suppose the assertion holds for $i=1,2,\ldots,m-1$. Note that $\theta_m=\theta_{m-1}+\theta ||\mQ(\theta_{m-1}||$
and by ${\cal E}^1$, $||\mQ(\theta_m)||\geq 2||\mQ(\theta_{m-1})||$. Therefore
\begin{align*}
{||\mQ(\theta_m)||\over \theta_m} &\geq {2||\mQ(\theta_{m-1})||\over \theta_{m-1}+\theta ||\mQ(\theta_{m-1})||} \\
&={2||\mQ(\theta_{m-1})||/\theta_{m-1}\over 1+\theta ||\mQ(\theta_{m-1})||/\theta_{m-1}}.
\end{align*}
But by the inductive assumption $||\mQ(\theta_{m-1})||/\theta_{m-1}\geq 1/\theta$. This immediately implies that the
expression above is also at least $1/\theta$, and the induction is completed. Now for every $t\geq \theta_0=\theta ||q||$ we find $\theta_m$
such that $\theta_{m-1}\leq t<\theta_m$. Using (\ref{eq:Qinfty}) we obtain
\begin{align*}
{||\mQ(t)||\over t}&\geq {||\mQ(\theta_{m-1})||\over 4\theta_m}\max\left({\gamma\over C},1\right) \\
&={1\over 4(\theta_{m-1}/||\mQ(\theta_{m-1})||)+4\theta }\max\left({\gamma\over C},1\right) \\
&\geq {\max\left({\gamma\over C},1\right)\over 8\theta}>0,
\end{align*}
where the last inequality follows since $||\mQ(\theta_{m-1})||/\theta_{m-1}\geq 1/\theta$. This
shows (\ref{eq:rate})

Now, suppose the event ${\cal E}_i$ fails to occur at some  $\theta_i$, and thus
${\cal E}^1$ does not occur. We then ``restart'' the process of
attempting to obtain an infinite sequence of points $\theta_i$ with the
properties outlined above. Let us call ${\cal E}^2$ the event that the
that sequence is obtained after restarting the process again as follows.
At that time at which ${\cal E}_i$ fails, we switch to any non-idling
non-preemptive scheduling policy. Applying Assumption C, with probability one
there exists a time $\tau_1$ for which $||\mQ(\tau_1)||\geq n_0$. Note that
it is possible that $\tau_1=\theta_i$.
We apply the policy ${\cal U}$ starting from time
$\tau_1$. Repeating the argument for ${\cal E}^1$, with probability greater than $\alpha$
we obtain a new
infinite sequence of time instances $\theta'_i$ such that
$||\mQ(\theta_{i+1}')||\geq 2||\mQ(\theta_i^{\prime})||$, i.e.\
${\cal E}^2$ occurs. If ${\cal E}^2$ does not occur we again restart the
process. Finally the probability of eventually obtaining a sequence of
points $\theta_i$ with the stated properties is given by
$\pr({\cal E}) \equiv \pr(\cup_{k=1}^{\infty}{\cal E}^k)$, where the ${\cal E}^k$
are defined as above in the natural way. Since the probability of each event
${\cal E}^k$ is bounded below by $\alpha$, and this lower bound on probability
does not depend on whether or not the other events occur, the
probability of ${\cal E}$ is one.

Finally, we show below that the event ${\cal E}$ implies (\ref{eq:rate0}), i.e.\
\be{eq:ratek}
\liminf_t{||\mQ(t)||\over t}\geq {\max\left({\gamma\over C},1\right)\over 8\max(1,{3\over \gamma})}>0.
\ee
Let $k$ be the smallest integer for which ${\cal E}_k$ occurs. Denote by $T_k$ the time corresponding
to the beginning of this event. Fix any state $(q,z_1,z_2)$ and $t_0>0$ and condition on $T_k=t_0, \mQ(t_0)=q$.
Applying (\ref{eq:rate}) to the event ${\cal E}_k$ we obtain
\be{eq:ratek1}
\liminf_t{||\mQ(t)||\over t_0+t}\geq {\max\left({\gamma\over C},1\right)\over 8\max(1,{3\over \gamma})}>0.
\ee
Since $t_0$ is fixed, the lower bound (\ref{eq:ratek}) holds as well. Integrating over the choices of $(q,z_1,z_2)$
and $t_0$ we complete the proof of the theorem. \qed

\section{Conclusions and Further Work}\label{section:conclusions}
The present work leaves many interesting questions open. The most immediate one
is whether the result connecting global weak stability and rate stability
holds for networks
with any number of stations. One way to prove this conjecture would be to
establish the Finite Decomposition
Property for fluid networks with more than two stations.
Of course the question of whether (strong) global stability of the fluid
model is equivalent to
positive Harris recurrence remains open even for networks with two
stations. There the difficulty
lies in begin able to analyze the dynamics of the stochastic network at the
critical regime $\rho^*=1$.
Finally, we mention that our assumption that interarrival and service times
are i.i.d. is
used to simplify the exposition and our result should hold for networks
with more
general primitives as long as the associated processes satisfy appropriate
large deviations bounds.

\section{Acknowledgements}
During the course of this research,
the second author was partially supported by NSF Grant DMI-0132038.
We would also like to thank several anonymous referees for helpful
technical and expositional comments.

\bibliographystyle{abbrv}
%% \bibliography{has}

\appendix

\section{Appendix}

\begin{proof}[Proof of Lemma~\ref{lemma:LDiid}] We begin by proving (\ref{eq:LDBtime}).
The proof of (\ref{eq:LDBrate}) is then derived using (\ref{eq:LDBtime}).
Our method uses the standard derivation of LD upper bounds on i.i.d.\
sequences.

\textbf{Part I.} Let us fix arbitrary $\epsilon>0$ and $\theta>0$. Then we
note that the following
hold for all $n \ge 1$:

\begin{eqnarray*}
\pr\left\{\sum_{1\leq i\leq n}Z_i\geq n\alpha+n\epsilon+z \; \bigg| \;
Z_1\geq z\right\} & =
& \pr\{e^{\theta\sum_{1\leq i\leq n}Z_i}\geq
e^{\theta(n\alpha+n\epsilon+z)}\; | \; Z_1\geq z\} \\
& \leq & {\mathbb{E}[e^{\theta(Z_1-z)}|Z_1\geq z](\mathbb{E}[e^{\theta
Z_2}])^{n-1}\over e^{n\theta(\alpha+\epsilon)}} \\
& \leq & {F(\theta)(\mathbb{E}[e^{\theta Z_2}])^{n}\over
e^{n\theta(\alpha+\epsilon)}},
\end{eqnarray*}
where we use $z\geq 0$ and $\mathbb{E}[e^{\theta Z_2}]\geq 1$. It
is a standard result in large deviations theory \cite{demzei98}
that $\mathbb{E}[e^{\theta Z_2}]/e^{\theta(\alpha+\epsilon)}\equiv
e^{-L(\epsilon)}<1$ for some value of $\theta=\theta(\epsilon) \in
[0,\theta_0]$ as long as $\mathbb{E}[e^{\theta Z_2}]$ is finite on
$[0,\theta_0]$. Thus our tail probability is at most
$F(\theta(\epsilon))e^{-L(\epsilon)n}$. We fix a suitable $\theta$ and take
$V=F(\theta)$.

We now prove a complimentary bound. Again fix arbitrary $\epsilon>0$ and
$\theta>0$.
\begin{eqnarray*}
\pr\left\{\sum_{1\leq i\leq n}Z_i\leq n\alpha-n\epsilon+z \; \bigg| \;
Z_1\geq z\right\} & =
& \pr\{e^{-\theta\sum_{1\leq i\leq n}Z_i}\geq
e^{-\theta n\alpha+\theta n\epsilon-\theta z}\; | \; Z_1\geq z\} \\
& \leq & {\mathbb{E}[e^{-\theta(Z_1-z)}|Z_1\geq z](\mathbb{E}[e^{-\theta
Z_2}])^{n-1}\over e^{-n\theta\alpha+n\theta\epsilon}} \\
& \leq & {(\mathbb{E}[e^{-\theta Z_2}])^{n-1}\over
e^{-(n-1)\theta(\alpha-\epsilon)}}e^{\theta\alpha-\theta\epsilon},
\end{eqnarray*}
where we use $\mathbb{E}[e^{-\theta(Z_1-z)}|Z_1\geq z]\leq 1$. Again we use
a standard result in large deviations theory \cite{demzei98}
stating that $\mathbb{E}[e^{-\theta
Z_2}]/e^{-(\theta\alpha-\theta\epsilon)}\equiv
e^{-L(\epsilon)}<1$ for some value of $\theta=\theta(\epsilon) \in
[0,\theta_0]$. We take $V=e^{L(\epsilon)+\theta\alpha-\theta\epsilon}$.
This proves (\ref{eq:LDBtime}).

\textbf{Part II.} We now prove (\ref{eq:LDBrate}). Consider a fixed,
but arbitrary $\epsilon > 0$.

We first obtain a bound which is valid for all $t \ge 1$.
Note that (\ref{eq:LDBtime}) easily implies
two one-sided versions of the LD inequality. Applying
one such one-sided version of (\ref{eq:LDBtime})
with $n = \lceil t/\alpha + \epsilon t \rceil \ge 1$ and $\tilde{\epsilon}=
\alpha^2 \epsilon / (1+\alpha \epsilon + \alpha) > 0$, then there exist
$\tilde{L}, V_1 > 0$ such that:
$$ \pr \left\{ \sum_{i=1}^{\lceil t/\alpha + \epsilon t \rceil}
Z_i -z \le \alpha \lceil t/\alpha + \epsilon t \rceil - \tilde{\epsilon}
\lceil t/\alpha + \epsilon t \rceil \; \bigg | \; Z_1 \ge z \right\} \le
V_1e^{-\tilde{L} \lceil t/\alpha + \epsilon t \rceil},$$
for all $t \ge 0$.
Next since $ \lceil t/\alpha + \epsilon t \rceil \ge t/\alpha + \epsilon t
$
and for $t \ge 1$, $\lceil t/\alpha + \epsilon t \rceil \le
t / \alpha + \epsilon t + t$ we have
$$ \pr \left\{ \sum_{i=1}^{\lceil t/\alpha + \epsilon t \rceil}
Z_i -z \le \alpha [t/\alpha + \epsilon t ] - \tilde{\epsilon}
[ t/\alpha + \epsilon t + t ] \; \bigg | \; Z_1 \ge z \right\} \le
V_1e^{-\tilde{L} \lceil t/\alpha + \epsilon t \rceil},$$
for all $t \ge 1$.
Multiplying through inside the probability yields:
$$ \pr \left\{ \sum_{i=1}^{\lceil t/\alpha + \epsilon t \rceil}
Z_i -z \le t + \alpha \epsilon t  - \alpha \epsilon t
\; \bigg | \;  Z_1 \ge z \right\} \le
V_1e^{-\tilde{L} \lceil t/\alpha + \epsilon t \rceil}.$$
Further simplification gives:
$$ \pr \left\{ \sum_{i=1}^{\lceil t/\alpha + \epsilon t \rceil}
Z_i  \le t + z
\; \bigg | \; Z_1 \ge z \right\} \le
V_1e^{-\tilde{L} \lceil t/\alpha + \epsilon t \rceil},$$
for all $t \ge 1$.
Using the duality relationship between a counting process and
its increments, the above implies:
\begin{equation} \label{eq:bandl}
\pr \left\{ N(t + z) \ge t/\alpha + \epsilon t \; | \; Z_1 \ge z \right\} \le
V_1e^{-\tilde{L} \lceil t/\alpha + \epsilon t \rceil}.
\end{equation}
Finally setting we have
$L_1 = \tilde{L} (1/\alpha + \epsilon)$,
$$   V_1e^{-\tilde{L} \lceil t/\alpha + \epsilon t \rceil} \le V_1 e^{-L_1
t }, $$
for all $t \ge 1$. Hence, we can rewrite (\ref{eq:bandl}) as
\begin{equation} \label{eq:l1}
\pr \left\{ N(t + z) \ge t/\alpha + \epsilon t \; | \; Z_1 \ge z \right\} \le
V_1 e^{-L_1 t},
\end{equation}
for all $t \ge 1$.

For $t < 1$ and any $V_2 >1$ note that $V_2 \cdot \mathrm{exp}(-L_2 t) \ge
1$ if
we set $L_2 \equiv \ln V_2 >0$. Hence,
\begin{equation} \label{eq:l2}
  \pr \left\{ N(t + z) \ge t/\alpha + \epsilon t \; | \; Z_1 \ge z \right\}
\le
V_2 e^{-L_2 t},
\end{equation}
holds trivially for all $t < 1$ with such a $V_2$ and $L_2$.
Finally, setting $V_3 \equiv \max\{V_1,V_2\}$ and
$L_3 \equiv \min \{L_1, L_2 \}$, then combining with (\ref{eq:l1}) and
(\ref{eq:l2})
we obtain:
$$ \pr \left\{ N(t + z) \ge t/\alpha + \epsilon t \; | \; Z_1 \ge z \right\}
\le
V_3 e^{-L_3 t},$$
for all $t \ge 0$.

This proves one side of the inequality in (\ref{eq:LDBrate}).
The other direction
is proved by an exactly analogous argument. The final result is then
obtained by combining the two directions, applying Boole's inequality and
again using appropriate
constants $B$ and $L$.

\qed
\end{proof}

\end{document}